\documentclass[onefignum,onetabnum,onealgnum,oneeqnum]{siamart190516}
\usepackage[papersize={8.5in,11in},margin=1.0in]{geometry}
\usepackage{amssymb}
\usepackage{amsmath}
\usepackage{enumitem}
\usepackage{multirow}
\usepackage{microtype}
\usepackage{bm}
\usepackage{algorithm}
\usepackage[noend]{algpseudocode}
\makeatletter
\algrenewcommand\ALG@beginalgorithmic{\footnotesize}
\renewcommand{\ALG@name}{\small Algorithm}
\makeatother
\setlist{leftmargin=1.5em}
\crefname{section}{\S\!}{\S\!}
\crefname{subsection}{\S\!}{\S\!}
\crefname{figure}{Fig.\!}{Figs.\!}

\newcommand{\R}{\mathbb R}
\newcommand{\diagdots}{\raisebox{-0.2em}{$\ddots$}}
\def\trianglesymbol{\raisebox{-0.1em}{\rotatebox{90}{\scalebox{0.95}[0.75]{$\blacktriangle$}}}\hspace{0.08em}}
\def\squaresymbol{\raisebox{-0.05em}{\rotatebox{45}{\scalebox{0.55}{$\blacksquare$}}}}

\newcommand{\circint}{\mathop{\mathpalette\docircint\relax}\!\int}
\newcommand{\docircint}[2]{%
  \ifx#1\displaystyle
    \displaycircint
  \else
    \normalcircint{#1}%
  \fi
}
\newcommand{\displaycircint}{\rotatebox{-11}{\scalebox{1.5}{$\diamond$}}\mkern-19mu}
\newcommand{\normalcircint}[1]{\smallerc{#1}\ifx#1\textstyle\mkern-11.1mu\else\mkern-10.6mu\fi}
\newcommand{\smallerc}[1]{\vcenter{\hbox{\rotatebox{-9}{\scalebox{\ifx#1\textstyle 1\else 0.75\fi}{$\diamond$}}}}}

\makeatletter
\newcommand{\dotDelta}{{\vphantom{\Delta}\mathpalette\d@tD@lta\relax}}
\newcommand{\d@tD@lta}[2]{\ooalign{\hidewidth$\m@th#1\mkern-1mu\cdot$\hidewidth\cr$\m@th#1\Delta$\cr}}
\makeatother

\makeatletter
\newcommand{\algmargin}{\the\ALG@thistlm}
\makeatother
\algnewcommand{\parState}[1]{\State\parbox[t]{\dimexpr\linewidth-\algmargin}{\strut #1\strut}}

\title{High-Order Quadrature on Multi-Component Domains\\Implicitly Defined by Multivariate Polynomials}

\author{Robert I.~Saye\thanks{Mathematics Group, Lawrence Berkeley National Laboratory, Berkeley, CA 94720 (\texttt{rsaye@lbl.gov})}}
\date{\today}

\begin{document}

\setlength{\textfloatsep}{1em}

\maketitle

\begin{abstract}
A high-order quadrature algorithm is presented for computing integrals over curved surfaces and volumes whose geometry is implicitly defined by the level sets of (one or more) multivariate polynomials. The algorithm recasts the implicitly defined geometry as the graph of an implicitly defined, multi-valued height function, and applies a dimension reduction approach needing only one-dimensional quadrature. In particular, we explore the use of Gauss-Legendre and tanh-sinh methods and demonstrate that the quadrature algorithm inherits their high-order convergence rates. Under the action of $h$-refinement with $q$ fixed, the quadrature schemes yield an order of accuracy of $2q$, where $q$ is the one-dimensional node count; numerical experiments demonstrate up to $22$nd order. Under the action of $q$-refinement with the geometry fixed, the convergence is approximately exponential, i.e., doubling $q$ approximately doubles the number of accurate digits of the computed integral. Complex geometry is automatically handled by the algorithm, including, e.g., multi-component domains, tunnels, and junctions arising from multiple polynomial level sets, as well as self-intersections, cusps, and other kinds of singularities. A variety of numerical experiments demonstrates the quadrature algorithm on two- and three-dimensional problems, including: randomly generated geometry involving multiple high-curvature pieces; challenging examples involving high degree singularities such as cusps; adaptation to simplex constraint cells in addition to hyperrectangular constraint cells; and boolean operations to compute integrals on overlapping domains.
\end{abstract}


\section{Introduction}

In this paper, we develop high-order quadrature methods for curved surfaces and volumes whose geometry is implicitly defined by the level sets of (one or more) multivariate polynomials. In particular, given a set of polynomial functions $\phi_i: \R^d \to \R$, $i = 1, \ldots, n$, we consider volume and surface integrals of the form
\begin{equation} \label{eq:intro} \int_{U \setminus \Gamma} f\,dx \quad \text{and} \quad \int_{U \cap \Gamma} g\,dS, \end{equation}
where the implicitly defined interface, $\Gamma := \bigcup_{i=1}^n \{ \phi_i = 0 \}$, is the union of each polynomial's zero level set (its ``interface''), and $U \subset \R^d$ is a given hyperrectangle (e.g., a rectangle in two dimensions, or rectangular box in three dimensions). Here, the integrands $f$ and $g$ need not be smooth across the entire hyperrectangle $U$; we instead assume that $f$ and $g$ are sufficiently smooth on each connected component of $U \setminus \Gamma$ and $U \cap \Gamma$, respectively. As such, our main problem formulation concerns a hyperrectangular domain which has been divided into one or more pieces by the zero level set of each polynomial, and we develop quadrature algorithms which can accurately compute integrals on each individual piece.

Integrals of the form \eqref{eq:intro} arise in a variety of applications involving implicitly defined geometry, including, e.g.: discretization methods for solving partial differential equations on domains with curved geometry, such as cut-cell, embedded boundary, immersed interface, and implicit mesh methods \cite{ImplicitMeshPartOne,ImplicitMeshPartTwo,JohanssonLarson,MoesDolbowBelytschko,ChengFries,BurmanClausHansboLarsonMassing}; %
remapping in arbitrary Lagrangian-Eulerian methods \cite{doi:10.1137/17M1116453}; multi-mesh methods, which employ overlapping meshes of different characteristics for different parts of a multi-physics problem \cite{JOHANSSON2019672}; and in the treatment of jump conditions and singular source terms \cite{flame,immersedboundary,Brackbill1992335,DiscacciatiQuarteroniQuinodoz}. %
In many such applications, the interfacial geometry is not known ahead of time and is instead computed dynamically, e.g., through the level set method \cite{OsheSeth,SethBook,OsherBook,LevelSetFluidReview}. It is often important to compute these integrals to a high degree of accuracy, e.g., to ensure consistency in a weak formulation or to sufficiently diminish the impact of variational errors. Further, the implicitly defined geometry can be complex, even in a single grid cell, having characteristics such as tunnels, junctions, corners, or cusps; in these cases, it would be ideal to directly treat the geometry at hand, without having to resort to heuristic subdivision/refinement algorithms. %

The quadrature algorithms developed in this work are designed with these aspects in mind and have the following features:
\begin{itemize}
\item The output is a quadrature scheme of the form
\[ \int_{U \setminus \Gamma} f \approx \sum_i w_i f(x_i) \quad \text{and} \quad \int_{U \cap \Gamma} g \approx \sum_j {\mathsf w}_j g({\mathsf x}_j), \]
where the weights $w_i$, ${\mathsf w}_j$ are strictly positive and the nodes  $x_i, {\mathsf x}_j$ are strictly inside their respective domains. Importantly, the quadrature schemes yield accurate methods for computing integrals on the individual connected components of these domains, e.g., if one polynomial $\phi$ implicitly defines $\Gamma$, then $\smash{\sum_{i | \phi(x_i) < 0} w_i f(x_i)}$ constitutes an accurate quadrature scheme for $\smash{\int_{U \cap \{ \phi < 0 \}} f}$, and there is no need to define or extend $f$ throughout all of $U$.

\item The algorithm employs two kinds of one-dimensional quadrature, Gauss-Legendre and tanh-sinh, and inherits their high-order convergence rates. In $d$ dimensions, based on a one-dimensional scheme of order $q$, the total number of quadrature nodes is ${\mathcal O}(q^d)$ in the volumetric case and ${\mathcal O}(q^{d-1})$ for surface integrals. Under the action of $h$-refinement with $q$ fixed, the schemes yield an order of accuracy approximately $2q$. Under the action of $q$-refinement with the geometry fixed, the convergence is approximately exponential in $q$, i.e., doubling $q$ approximately doubles the number of accurate digits. 

\item Complex domain topology is automatically detected by the algorithm, including, e.g., multi-component domains, tunnels, and intersections arising from multiple polynomial level sets, as well as corners, cusps, and other kinds of singularities. Instead of relying on an ad hoc subdivision/refinement process, the methods divide the integral domain into a finite number of (implicitly defined) pieces whose boundaries precisely align with the interfacial geometry.

\item The overall approach is dimension independent; in particular, the same approach can be used in two- and three-dimensional applications. %
\end{itemize}
A wide variety of accompanying numerical experiments assess the effectiveness of the quadrature algorithm. Examples include convergence studies under $h$-refinement and $q$-refinement; performance on randomly generated polynomials having non-trivial interfacial topology; ability to handle singularities such as cusps and self-intersections; adaptation to simplices instead of hyperrectangles; and junctions arising from overlapping implicitly defined domains. 

An outline is as follows. In the next section, we briefly review some typical strategies for computing integrals of implicitly defined domains, and then provide some context motivating the design of the algorithms in this work. Section \cref{sec:driver} presents the main algorithm, starting with an overview and then following with a discussion of individual driving components. A broad range of numerical tests is presented in \cref{sec:results}, and concluding remarks with a brief discussion are given in \cref{sec:conclusion}. Meanwhile, the supplementary material accompanying this paper (available online) includes additional numerical experiments and algorithms.

\section{Background}

A wide variety of methods have been developed to compute integrals on implicitly defined domains. General approaches typically involve one of the following strategies; in each case, we indicate some of the associated challenges.
\begin{itemize}
\item \textit{Methods which reconstruct an approximation of the interface.} For example, one approach is to build a piecewise linear approximation of the interface %
and then apply standard quadrature algorithms on the resulting simplices or polyhedra, %
see, e.g., \cite{MinGibou2007,HoldychNobleSecor,Ventura}. Applied to complex interface geometry, this approach can have limited effectiveness due to a relatively low approximation power; often, subdivision techniques are required to locally refine the geometry approximation, and often without a good return on the increased computational cost. %
Higher-order explicit reconstruction schemes are possible, see, e.g., \cite{doi:10.1137/17M1126370,ENGVALL201783,ChengFries,10.1002/nme.5121,FRIES2017759,BOCHKOV20191156}; while effective in some cases, the implementation of these methods quickly becomes more subtle as the approximation order is increased, in part due to a highly sensitive relationship between interface topology and Jacobian calculations. %

\item \textit{Perturbation and correction methods.} These methods begin with a primary interface approximation, typically piecewise linear, and then incorporate a correction term %
representing the effect of locally projecting the approximation to the true curved interface geometry,
see, e.g., \cite{LEHRENFELD2016716,SCHOLZ2019112577,scholz2020high,THIAGARAJAN20141682,doi:10.1098/rspa.2016.0401}. In some cases, these methods offer a good compromise between the ease of piecewise linear reconstruction and the accuracy of a mildly higher-order interface representation. However, due to the reliance on the initial surface approximation, this kind of approach can breakdown on high curvature geometry or multiple surfaces in the same grid cell.

\item \textit{Application of smoothed Dirac delta and Heaviside functions.} These methods avoid an explicit reconstruction of the interface; instead, a sampling of the given level set function, typically on Cartesian grids, is used to approximate, e.g., $\int_\Gamma f = h^d \sum_i f(x_i) \delta_h(\phi(x_i)) |\nabla_h \phi(x_i)|$, where $\delta_h$ is a grid-dependent smoothed Dirac delta function and $x_i$ are the grid points. These approaches typically rely on a cancellation of errors in the summation over regularly spaced grid points and it can be subtle to develop convergent schemes \cite{TornbergEngquist,EngquistTornbergTsai,Smereka,Towers,ZahediTornberg,MinGibou2008}, especially higher-order ones \cite{Wen2008,Wen2009,Wen2010}; see also related methods built through application of the co-area formula \cite{doi:10.1137/16M1102227,KublikTsai,10.1115/1.4039639,10.1115/1.4047355}. %

\item \textit{Mechanisms making use of integration by parts.} Quite distinct to the above approaches, another possibility is to use integration by parts (e.g., via the divergence theorem or the generalized Stokes theorem) to recast an implicitly defined volume/surface integral on the interior of a cell $U$ as a surface integral on the boundary of $U$. In general, the domain of the resulting surface integral is itself implicitly defined, and so this process may continue recursively, removing one dimension at a time. Under certain assumptions, these methods can yield accurate and effective results, see, e.g., moment-based methods \cite{MullerKummerOberlack,10.2140/camcos.2015.10.83}, their application to planar regions with piecewise parameterised boundaries \cite{SOMMARIVA2009886,GUNDERMAN2021102944}, and trimmed parametric surfaces \cite{gunderman2021highaccuracy}. However, as discussed in \cite{HighorderImplicitQuad}, methods based on integration by parts may not always work, even for highly-resolved geometry---there is an implicit assumption that the boundary integrals on $\partial U$ can detect enough of the internal geometry; in some cases, these methods require subdivision to prevent parts of the geometry from being missed.

\item \textit{Recasting as the graph of a height function.} Another possibility is to recast the implicitly defined geometry as the graph of an implicitly defined height function; using this interpretation, a quadrature scheme can then be formulated through the corresponding splitting of coordinates, i.e., through two nested integrals, one in the vertical/height direction, and the other in all remaining tangential directions. In the author's prior work \cite{HighorderImplicitQuad,algoim}, this approach was used to build very high-order quadrature schemes, e.g., 20th order; see also its extension to simplices by Cui \textit{et al} \cite{10.1145/3372144}. These methods can be fast, requiring only one-dimensional root finding, and inherit the high-order accuracy of the underlying one-dimensional Gaussian quadrature schemes. However, the overall approach relies on finding a suitable coordinate axis such that the interface is the graph of a well-defined, sufficiently smooth height function, e.g., without any corners, cusps, or vertical tangents. In \cite{HighorderImplicitQuad,algoim}, interval arithmetic is used to subdivide the hyperrectangle $U$, if necessary, until the level set function is provably monotone in a particular height direction $k$, thereby guaranteeing (by the implicit function theorem) a well-defined height function. This process is effective for a wide kind of geometry, but can fail in some cases, e.g., in the presence of singularities, corners, or junctions. In these situations, an excessive amount of subdivision can result, a process that may also need to be halted, potentially leaving behind small voids. 
\end{itemize}

\noindent Putting aside the challenging cases of complex interface geometry, the above approaches can typically be applied to a broad class of sufficiently smooth level set functions. In the more restrictive case of polynomial level set functions, there is more structure that can be exploited, e.g., to effectively treat singularities. However, there appears to be relatively few works in the literature devoted specifically to building quadrature schemes for domains implicitly defined by polynomials. In some cases, these efforts have concentrated solely on computing the volume of semi-algebraic sets. For example, in Lairez \textit{et al} \cite{10.1145/3326229.3326262}, and with a focus on arbitrary-precision calculations in algebraic geometry, the authors develop a method to compute volumes of compact semi-algebraic sets; the methods apply a dimension-reduction process to reformulate the volume computation as a solution of a high-order Picard-Fuchs differential equation, which in turn relies on high-precision ODE solvers. As another example, Henrion \textit{et al} \cite{doi:10.1137/080730287} show how the volumes and moments of basic semi-algebraic sets can be approximated using a hierarchy of semi-definite programming problems, with guaranteed lower and upper bounds on these estimates, though the accuracy of these methods appear limited. A third example is of Hrivn\'ak \textit{et al} \cite{Hrivn_k_2016}, which build high-accuracy quadrature schemes using very few points, but only for integrals involving a specific kernel and for specific domains.

Our motivation in this work is to address some of the above-mentioned challenges, particularly in regard to handling complex interface geometry involving high degrees of curvature, non-trivial topology with multiple components, and characteristics such as corners, intersections, and cusps. To do so, we focus on implicitly defined domains given by multivariate polynomials, and then exploit the additional structure this representation affords. In particular, the quadrature algorithms developed in this work have the following design objectives in mind:
\begin{enumerate}[label=(\roman*),leftmargin=3em]
\item Conceptually independent of the number of dimensions; e.g., the algorithms can be used for both two- and three-dimensional problems.
\item As close as possible to a ``black box'' algorithm; e.g., for use with dynamically-evolving interfacial geometry.
\item \label{nosub} Avoid the use of subdivision; in particular, subdivision is permitted only when it can be proven that a finite (usually small) number of subregions is needed.
\item \label{singular} Effectively handle cases of junctions, corners, cusps, etc., ideally with near-exponential convergence in the number of quadrature points.
\end{enumerate}
To achieve these objectives, we have adopted here a dimension reduction process in which the implicitly defined geometry is reinterpreted as the graph of an implicitly defined, potentially multi-valued height function. This approach is similar to the one used in the author's prior work \cite{HighorderImplicitQuad}, however the methods therein required a single-valued, well-defined and smooth height function. The objectives above, in particular \ref{nosub} and \ref{singular}, preclude this possibility and instead necessitate multi-valued height functions having additional kinds of characteristics, which in turn require additional mechanisms to effectively handle, as discussed next. %

\section{Quadrature on Domains Implicitly Defined by Multivariate Polynomials}
\label{sec:driver}

\subsection{Main formulation}
\label{sec:mainformulation}
 
As motivated in the prior section, our essential approach is to apply a dimension reduction technique wherein the implicitly defined interface $\Gamma$ is recast as the graph of an implicitly defined multi-valued ``height'' function. Once a suitable coordinate axis representing the vertical height direction is chosen, a quadrature scheme is built by performing: (a) an integration in the height direction along with (b) an integration in the tangential direction(s). Part (a) applies one-dimensional quadrature schemes to compute vertical line integrals, subdivided into pieces according to the location of the interface, while part (b) is itself recast as an integral over implicitly defined regions, in one fewer dimensions, whose associated integrand evaluates part (a). The latter integral is called the ``base integral'' and its domain shall be implicitly defined by a collection of multivariate polynomials; as such, the quadrature algorithm can execute itself and proceed recursively, reducing the number of dimensions one at a time. In order to obtain high quadrature accuracy, it is important to appropriately identify and separate the different regions of behaviour of the multi-valued height function. Changes in behaviour typically correspond to wherever the interface crosses the upper or lower faces of the constraint domain $U$, wherever two surfaces cross each other, or at branch points for which two surfaces join together at a vertical tangent. These aspects will be discussed shortly, but first we establish some preliminaries. %

Throughout this section we take, for ease of exposition, a reference domain equal to the $d$-dimensional unit cube; adaptation to any hyperrectangular domain $U \subset \R^d$ is a straightforward modification. For the time being, we will also focus attention on the case of volumetric integrals; surface integrals will be treated at the end of the section in \cref{sec:surface}. As such, suppose we are given one or more polynomial functions $\phi_i : (0,1)^d \to \R$, $i = 1, \ldots, n$, which together implicitly define an interface $\Gamma := \bigcup_{i=1}^n \{\phi_i = 0\}$. Suppose also that we are given an integrand $f : (0,1)^d \to \R$, which is permitted to exhibit discontinuities across the interface, but is otherwise assumed to be sufficiently smooth. For the moment, let us also suppose we have identified a suitable coordinate axis representing the height direction, $e_k$, say.\footnote{$e_i$ denotes the standard basis vector in the direction of the $i$th coordinate.} Our task is to build a quadrature scheme to integrate $f$ over the volumetric region on either side of $\Gamma$, constructed by performing an outer integral (over $(0,1)^{d-1}$) of an inner 1D integral in the vertical direction $e_k$.

Denote by $x = (x_1, \ldots, x_{k-1}, x_{k+1}, \ldots, x_d) \in (0,1)^{d-1}$ a point in the $(d-1)$--dimensional cube, and denote by $x + y e_k$ the point $(x_1, \ldots, x_{k-1}, y, x_{k+1}, \ldots, x_d)$, i.e., the lifting of $x$ in the direction $e_k$ by the amount $y$.  We denote the integrand of the base integral by ${\mathcal F} : (0,1)^{d-1} \to \mathbb R$. In particular, given a base point $x$, evaluation of ${\mathcal F}$ at $x$ performs an integral over a vertical line segment $\{ x + y e_k : 0 \leq y \leq 1 \}$; since our goal is to build high order quadrature schemes which can separately handle either side of $\Gamma$, it follows that the inner integral should be divided into pieces, demarcated by interfacial crossing points. This can succinctly be described by asking the inner integral to build a quadrature scheme which integrates $f$ on the connected components of the set $(0,1) \setminus \Gamma(x)$, where $\Gamma(x) := \bigcup_{i=1}^n \{y : 0 \leq y \leq 1 \text{ and } \phi_i(x + y e_k) = 0 \}$. Here, it is convenient to introduce a custom notation: denote by an integral symbol with overset $\rotatebox{-11}{\scalebox{1.5}{$\diamond$}}$ as summing over the connected components of the corresponding domain:
\begin{equation} \label{eq:cint} \circint_\Omega := \sum_j \int_{\Omega_j} \text{ where $\Omega_j$ enumerates the connected components of $\Omega$.} \end{equation} 
Using this notation, the base integral's integrand ${\mathcal F} : (0,1)^{d-1} \to \mathbb R$ is given by
\begin{equation} \label{eq:innerintegral} {\mathcal F}(x) := \circint_{(0,1) \setminus \Gamma(x)} f(x + y e_k)\,dy \end{equation}
and part of the algorithm's role is to produce a quadrature scheme which approximately evaluates ${\mathcal F}(x)$ for a given $x$. %

\begin{figure}%
\centering\includegraphics{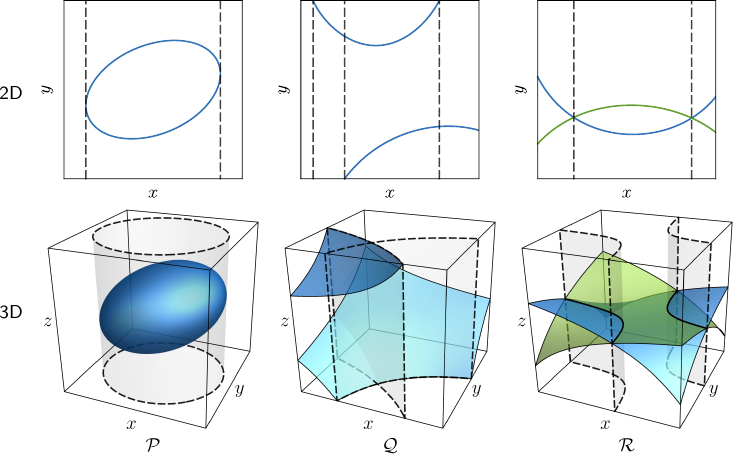}%
\caption{Depictions of the changes in behaviour of $\Gamma$ according to its representation as the graph of a multi-valued height function. In each case, the height direction is ``up'' (i.e., $\hat y$ in 2D and $\hat z$ in 3D), the blue and green lines/surfaces illustrate the interface, and the dashed lines/surfaces represent the vertical extrusion of the sets $\mathcal P$, $\mathcal Q$ and $\mathcal R$, from the base cube $(0,1)^{d-1}$ into the whole cube $(0,1)^d$.  \textit{(left column)} $\mathcal P$ denotes the set of base points above which the interface branches, typically occurring at the vertical tangents of rounded surfaces. \textit{(middle column)} $\mathcal Q$ denotes the set of base points above which the interface crosses the upper or lower faces of $(0,1)^d$, i.e., wherever the multi-valued height function has one or more values equal to $0$ or $1$. \textit{(right column)} $\mathcal R$ denotes the set of base points above which two different interfaces cross.}
\label{fig:examples}
\end{figure}

Note that the smoothness of ${\mathcal F}(x)$ depends on changes in the connected component behaviour of $\Gamma(x)$, which are precisely the changes in the behaviour of $\Gamma$ according to its representation as the graph of a multi-valued height function. To illustrate, \cref{fig:examples} demonstrates some typical scenarios in 2D and 3D. In each case, the coordinate axis representing the vertical height direction is ``up'', i.e., $\hat y$ in 2D and $\hat z$ in 3D, and the points of interest are subsets of the base cube $(0,1)^{d-1}$. Points labelled $\mathcal Q$ indicate where the interface $\Gamma$ crosses the upper or lower faces, and correspond to where one or more values of the height function equals $0$ or $1$. Points labelled $\mathcal R$ indicate the base points above which two different interfaces cross. Points labelled $\mathcal P$ indicate base points above which the height function branches; typically, branching occurs at vertical tangents of rounded surfaces such as circles and spheres. 

Combined, these points mark the boundaries of the piecewise-smooth behaviour of ${\mathcal F}(x)$. To attain high-order accuracy in the overall quadrature scheme, it follows that the base integral should be prevented from ``jumping across'' these boundaries. We achieve this through a division of the base integral's domain into one or more pieces, such that the base integral produces an accurate quadrature scheme for each piece. As before, this can be succinctly formulated through the use of connected components of an implicitly defined domain; specifically, we begin with the base cube $(0,1)^{d-1}$ and then subtract out an implicit description of the points $\mathcal P \cup \mathcal Q \cup \mathcal R$ discussed above; in essence, the points $\mathcal P \cup \mathcal Q \cup \mathcal R$ form the implicitly defined interface of the $(d-1)$-dimensional base integral.

Implicitly characterising points in $\mathcal Q$, i.e., changes in the smoothness of ${\mathcal F}$ arising from interfacial crossings at the upper and lower faces, is straightforward: this can be accomplished by excluding the sets $\{x: \phi_i(x + 0 e_k) = 0\}$ and $\{x: \phi_i(x + 1 e_k) = 0\}$ from the base integral's domain.

An implicit characterisation of $\mathcal P$ and $\mathcal R$ requires additional tools, namely \textit{discriminant} and \textit{resultant} methods. A key tool in algebraic geometry, these methods take in as input one or two polynomials and output a new polynomial whose roots pinpoint the location of the shared roots of the input (and its derivative, in the case of discriminants). More details on their use will be given in \cref{sec:resultants}; for now, it suffices to state their application yields multivariate polynomials, in one fewer dimensions, whose zero level sets contain the set of points $\mathcal P$ and $\mathcal R$. These polynomials, in conjunction with the upper- and lower-face restrictions of $\phi_i$ used to characterise $\mathcal Q$, complete the partitioning of the base integral's domain into one or more connected components.

The output of the above analysis is a set of polynomials $\psi_j : (0,1)^{d-1} \to \R$ describing the domain of the base integral in the form $(0,1)^{d-1} \setminus \bigcup_j \{\psi_j = 0\}$. %
Using the notation of \eqref{eq:cint}, the base integral then takes the form
\[ \circint_{(0,1)^{d-1} \setminus \bigcup_j \{\psi_j = 0\}} {\mathcal F}(x)\, dx. \]
Note that this is of the same form as we started with, but in one fewer dimensions: the task of the base integral is to build a quadrature scheme for an implicitly defined domain, whose interface is the union of a collection of polynomial isosurfaces, and whose integrand is, by construction, smooth on either side of that interface. Therefore, we can apply the quadrature algorithm recursively, eliminating one coordinate direction at a time, until we reach the conceptually simple one-dimensional base case.

We are nearly in a position to give a high-level description of the quadrature algorithm, however some aspects need additional motivation: 
\begin{itemize}
\item \textit{Quadrature methods.} By construction, we need only apply one-dimensional quadrature schemes, but there is a choice in the specific kind of scheme to employ. Indeed, it can be worthwhile to apply different kinds of 1D quadrature in different parts of the dimension-reduction process. In this work we have mainly explored the use of Gauss-Legendre and tanh-sinh quadrature, choosing between the two depending on the smoothness of the height functions involved. This choice is one of the main topics of discussion in this paper, and some guidelines will be given.

\item \textit{Polynomial basis.} Numerically, a basis for multivariate polynomials must be chosen. Given the dimension-reduction process and our choice of working with hyperrectangular domains, a natural candidate is a tensor-product basis of polynomials. In particular, we have made extensive use of the Bernstein basis, owing to its favourable numerical properties such as well-conditioned root finding and accurate range/bounds evaluation; more comments are given in \cref{sec:bernstein}.

\item \textit{Polynomial elimination and masking.} It can be highly beneficial to detect when polynomials can be eliminated in the dimension reduction process. For example, if we can determine that $\Gamma$ is empty, then a simple tensor-product quadrature scheme would suffice. As another example, if we can detect that no branching points exist (e.g., points labelled $\mathcal P$ in \cref{fig:examples}), then two benefits arise: (i) we can avoid calculating discriminants, and (ii) conclude that the associated height function has no vertical tangents, thereby helping to decide which quadrature scheme to use. In this work, these kinds of checks are performed with efficient ``mask'' operations which, in turn, exploit useful properties of the Bernstein basis. Roughly speaking, a mask divides the cube $(0,1)^d$ into a regular grid of $\mathcal M \times \cdots \times \mathcal M$ of subcells (typically $\mathcal M = 4$ or $8$); on each subcell, the mask has a binary value of $0$ and $1$ and indicates whether its accompanying polynomial is provably nonzero on the subcell or if its roots can be ignored. More details on masks and their use are given in \cref{sec:masks}.

\item \textit{Choice of height direction/elimination axis.} Prior to applying the dimension-reduction process, one must make an appropriate choice for the height direction, i.e., the coordinate axis to be eliminated, and in most applications, this choice will need to be made algorithmically. A variety of aspects could impact this decision, including: which directions yield height functions with the smallest amount of curvature (and thus highest quadrature accuracy); which directions yield branching points/vertical tangents (which should generally be avoided, if possible); and which directions lead to the fewest number of polynomials in the description of the base integral's implicitly defined domain (which can benefit efficiency). Some of these aspects can be assessed with the masking operations. Further remarks are given in \cref{sec:height}.
\end{itemize}

\begin{algorithm}[t]
\caption{\small Compute $\mathcal{I} (f; \{(\phi_i, m_i)\}_{i=1}^n, q)$. Given a list of multivariate polynomials $\phi_i : (0,1)^d \to \R$ and associated masks $m_i : {\mathbb N}_{\mathcal M}^d \to \{0,1\}$, a function $f : (0,1)^d \to \R$, and the order of quadrature $q$, approximate the volume integral $\circint_{(0,1)^d \setminus \Gamma} f(x)\,dx$, where $\Gamma := \bigcup_i \{\phi_i = 0\}$.}
\label{algo:main}
\begin{algorithmic}[1]

	\If{$d = 1$}
		\State Execute one-dimensional base case using the empty set as input base point:
		\State \textbf{return} ${\mathcal F}_\mathsf{eval}(\varnothing; 1, f, \{(\phi_i, m_i)\}_{i=1}^n, q)$.
	\EndIf
	
	\State Choose a height direction/elimination axis $k \in \{1, \ldots, d\}$ (see \cref{sec:height}). \label{algo:mainchoice}
	\State Define $\Psi$ to be a list of polynomial and mask pairs, initially empty.
	
	\Statex \textit{\hspace{-1em}Examine the lower and upper face restrictions of input polynomials:}	
	\For{each $i$}
		\For{$y \in \{0, 1\}$}
			\State Define $\psi : (0,1)^{d-1} \to \R$ by the face restriction of $\phi_i$, i.e., $\psi(x) = \phi_i(x + y e_k)$.
			\State Compute $m := \verb|nonzeroMask|(\psi, \verb|faceRestriction|(m_i, k, y))$.
			\If{$m$ is nonempty}
				\State Add $(\psi,m)$ to the list $\Psi$.
			\EndIf
		\EndFor		
	\EndFor
	
	\Statex \textit{\hspace{-1em}Examine the pseudo-discriminants of input polynomials:}
	\For{each $i$}
		\State Compute the polynomial $\partial_k \phi_i$.
		\State Compute $m := \verb|intersectionMask|\bigl((\phi_i, m_i), (\partial_k \phi_i, m_i)\bigr)$.
		\If{$m$ is nonempty}
			\State Compute the pseudo-discriminant $\psi := \dotDelta(\phi_i; k)$.
			\State Add $(\psi, \verb|collapseMask|(m, k))$ to the list $\Psi$.
		\EndIf
	\EndFor
	
	\Statex \textit{\hspace{-1em}Examine the pairwise resultants of input polynomials:}
	\For{each pair $i$ and $j$ with $i \neq j$, ignoring order}
		\State Compute $m := \verb|intersectionMask|\bigl((\phi_i, m_i), (\phi_j, m_j)\bigr)$. 
		\If{$m$ is nonempty}
			\State Compute the resultant $\psi := R(\phi_i, \phi_j; k)$.
			\State Add $(\psi, \verb|collapseMask|(m, k))$ to the list $\Psi$.
		\EndIf
	\EndFor
	
	\Statex \textit{\hspace{-1em}Apply dimension reduction:}
	\State Define a new integrand ${\mathcal F}: (0,1)^{d-1} \to \R$ by ${\mathcal F}(x) := {\mathcal F}_\mathsf{eval}(x; k, f, \{(\phi_i, m_i)\}_{i=1}^n, q)$.	\label{algo:mainintegranddef}
	\State \textbf{return} $\mathcal{I} ({\mathcal F}; \Psi, q)$.
	
	\end{algorithmic}
\end{algorithm}

\begin{algorithm}[t]
\caption{\small Compute $\mathcal{F}_\mathsf{eval}(x; k, f, \{(\phi_i, m_i)\}_{i=1}^n, q)$. Given a point $x \in \R^{d-1}$, a height direction $k$, a function $f : \R^d \to \R$, a list of multivariate polynomials $\phi_i : (0,1)^d \to \R$ and associated masks $m_i : {\mathbb N}_{\mathcal M}^d \to \{0,1\}$, and the order of quadrature $q$, approximate the one-dimensional integral $\circint_{(0,1) \setminus \Gamma(x)} f(x + y e_k)\,dy$, where $\Gamma(x) := \bigcup_i \{y : 0 \leq y \leq 1 \text{ and } \phi_i(x + ye_k) = 0\}$.}
\label{algo:integrand}
\begin{algorithmic}[1]
	\State Define $\mathcal R$ to be a list, initialised by $\mathcal R := \{0, 1\}$.
	\For{each $i$}
		\State Define $\psi: (0,1) \to \R$ by the vertical line restriction of $\phi_i$ above $x$ in the direction $k$, i.e., $\psi(y) := \phi_i(x + y e_k)$.
		\State Compute the set of real roots of $\psi$ in the interval $(0,1)$.
		\For{each such root $y$}
			\State Let $j \in {\mathbb N}_{\mathcal M}^d$ denote the mask subcell containing the point $x + y e_k$.
			\If{$m_i(j) = 1$}
				\State Add $y$ to the list $\mathcal R$.
			\EndIf
		\EndFor
	\EndFor
	\State Sort $\mathcal R = \{r_1, \ldots, r_\ell\}$ into ascending order such that $0 = r_1 < r_2 \leq \cdots \leq r_{\ell - 1} < r_\ell = 1$.
	\State Set $I = 0$.
	\For{$j = 1$ \textbf{to} $\ell - 1$}
		\If{$[r_j, r_{j+1}]$ is non-degenerate}
			\parState{On this interval, apply the chosen 1D quadrature scheme of order $q$ (e.g., Gauss-Legendre or tanh-sinh); let $x_{q,i}$ and $w_{q,i}$ denote the quadrature nodes and weights, respectively, relative to the interval $(r_j, r_{j+1})$:}
			\State Update $I := I + \sum_{i=1}^q w_{q,i} f(x + x_{q,i} e_k)$.
		\EndIf
	\EndFor
	\State \textbf{return} $I$.
	\end{algorithmic}
\end{algorithm}

With these preliminaries, a high-level description of the overall quadrature scheme is given in \cref{algo:main} and \cref{algo:integrand}. Various core operations in these algorithms have yet to be discussed in detail; we shall do so in the forthcoming sections. Here, a few remarks are in order:
\begin{itemize}
\item \cref{algo:main} is the main engine and follows the dimension-reduction approach motivated above: it chooses a height direction/elimination axis $k$, formulates an implicitly-defined domain for the base integral, and creates a new integrand ${\mathcal F}$ for the base integral. Meanwhile, the role of \cref{algo:integrand} is to evaluate the constructed integrands. %
\item The dimension reduction process is started by executing \cref{algo:main} on the user-provided multivariate polynomials (often just one of them). As coordinate axes are eliminated one by one, eventually the base case in $d = 1$ dimensions is reached; the base case executes \cref{algo:integrand} which builds a 1D quadrature scheme and evaluates the associated integrand at these points; in turn, a new 1D quadrature scheme is executed, and this process continues recursively, adding one dimension at a time. By construction, the result is a tensor-product quadrature scheme warping to the geometry found along the way. Although these algorithms assumed that a top-level integrand $f : (0,1)^d \to \R$ was provided by the user, they can easily be modified to record the final quadrature scheme as a collection of points $x_i \in (0,1)^d$ with associated weights $w_i \in \R^+$. %
\item From an object-oriented point of view, one may interpret the first algorithm as constructing integrand function objects, i.e., functors, whose evaluation is defined by the second algorithm. In fact, many of the operations in \cref{algo:main} are independent of the evaluation of the quadrature scheme; as such, one can cache the results and execute with different quadrature orders $q$, different 1D quadrature schemes, etc.
\end{itemize}

\cref{algo:main,algo:integrand} represent the main driver routines for the quadrature algorithms developed in this work. In the remainder of this section, we discuss in more detail some of the tools and algorithms needed to fully realise their implementation. A summary of the development is then given at the conclusion of this section, prior to the presentation of an extensive range of test cases in \cref{sec:results}.

\subsection{Bernstein basis}
\label{sec:bernstein}

Numerically, a basis must be chosen for the multivariate polynomials being manipulated by \cref{algo:main,algo:integrand}. We have adopted here a Bernstein basis, owing to its many favourable properties:
\begin{itemize}
\item Computing roots of polynomials in the Bernstein basis is generally well-conditioned. In fact, for simple roots in $[0,1]$ of an arbitrary polynomial, the root condition number is smaller in the Bernstein basis than in any other non-negative polynomial basis, including the power/monomial basis \cite{10.1090/S0025-5718-96-00759-4,FAROUKI19881}. %
\item One of the most accurate (and straightforward) methods for range evaluation of a polynomial is through the Bernstein basis \cite{garloff1993bernstein,Rajyaguru2017,garloff1993bernstein,lin1995methods,beska1989convexity}. These range evaluation methods employ the convex hull property of Bernstein polynomials and guarantee that a polynomial's value, at any point its the rectangular reference domain, is no smaller (or larger) than its minimum (or maximum) coefficient; in fact, these bounds become monotonically more accurate under the (stable) operations of Bernstein degree elevation and subdivision \cite{garloff1993bernstein,lin1995methods}.
\item Simple and stable algorithms exist to manipulate polynomials in Bernstein form, including evaluation, differentiation, degree elevation, and subdivision, among a number of other useful operations; see, e.g., \cite{FAROUKI19881}.
\end{itemize}
For the dimension reduction algorithm considered in this work, it is natural to employ a tensor-product Bernstein basis. In $d$ dimensions, let $(n_1,\ldots,n_d)$ denote the degree of a Bernstein polynomial $\phi$; then $\phi$ takes the form
\begin{equation} \label{eq:tensorprod} \phi(x_1,\ldots,x_d) = \sum_{i_1 = 0}^{n_1} \cdots \sum_{i_d = 0}^{n_d} c_{i_1,\ldots,i_d} b_{i_1}^{n_1} (x_1) \cdots  b_{i_d}^{n_d} (x_d) \end{equation}
where $b_\ell^n(x) = \binom{n}{\ell} (1 - x)^{n-\ell} x^\ell$, $\ell = 0, \ldots, n$, are the Bernstein basis functions of degree $n$ and $c_i = c_{i_1,\ldots,i_d}$ are the Bernstein coefficients of $\phi$. This formulation represents $\phi$ in coordinates relative to the standard unit domain $[0,1]^d$; subdivision refers to transforming the coefficients relative to a subset hyperrectangle, and can be stably computed using the \textit{de Casteljau} algorithm. This operation, along with many other standard operations in Bernstein form (such as differentiation, degree elevation, multiplication, etc.) are usually given for one-dimensional Bernstein polynomials, see, e.g., \cite{FAROUKI19881}, but it is straightforward to extend them to multivariate tensor-product polynomials of mixed degrees. %

\subsection{Masking}
\label{sec:masks}

A number of operations in \cref{algo:main} benefit from a cheap and robust method to isolate zeros of a polynomial or to detect when two polynomials have non-intersecting zero isosurfaces; for example, if it can be shown that all of the given polynomials $\phi_i$ are nonzero, then $\Gamma = \bigcup_i \{\phi_i = 0\}$ is empty and a simple tensor-product quadrature scheme for $(0,1)^d$ would suffice. However, our main use for these root isolation methods are to: (i) detect when $\Gamma$ is the graph of a multi-valued height function devoid of branching (to inform the type of 1D quadrature to use); (ii) avoid the calculation of a resultant or discriminant (which is potentially expensive in comparison to other aspects of the algorithm); and (iii) mitigate unnecessary subdivision caused by unwanted roots of resultant and discriminant calculations. These aspects will motivated in subsequent sections. %

A \textit{mask} $m : {\mathbb N}_{\mathcal M}^d \to \{0,1\}$ is a binary-valued function defined on a $\mathcal M \times \cdots \times \mathcal M$ grid of subcells $\Omega_i$, where the $i$th subcell of $[0,1]^d$ is given by $\Omega_i = \prod_{\ell=1}^d [\tfrac{1}{\mathcal M} (i_\ell - 1),\tfrac{1}{\mathcal M} i_\ell]$.\footnote{${\mathbb N}_{\mathcal M} := \{1, \ldots, {\mathcal M}\}$ denotes the first $\mathcal M$ positive integers.} Roughly speaking, a mask is built from one or two polynomials and helps to isolate the interfacial geometry of those polynomial(s)---a subcell value of ``$0$'' means the subcell can be (provably) ignored, while a value of ``$1$'' indicates that the subcell likely contains interfacial geometry. More concretely, the binary values of a mask $m$ are related to its associated polynomial $\phi$ in the following way: $m(i) = 0$ if it can be proven that $\phi$ is nonzero on $\Omega_i \pm \epsilon$; $m(i) = 0$ if it can be proven that any zeros of $\phi$ in $\Omega_i$ can be ignored by vertical line integrals (see \cref{sec:resultants}); and $m(i) = 1$ otherwise. Here, $\Omega_j \pm \epsilon$ denotes a slightly expanded version of $\Omega_j$, overlapping with neighbouring subcells by a small amount $\epsilon$: this overlap is used to eliminate edge/border cases that could be affected by roundoff error in floating point arithmetic.

These masks are (cheaply) computed using convex hull properties of Bernstein polynomials and relate to their effectiveness in accurately computing polynomial ranges/bounds cited earlier in \cref{sec:bernstein}. The design of these algorithms is deferred to the supplementary material; here, we summarise the essential ideas pertinent to their use in \cref{algo:main,algo:integrand} and the following sections:
\begin{itemize}
\item Each polynomial in the dimension reduction process, whether it be upper- or lower-face restrictions, discriminants, or resultants, has a mask associated with it. Typically, the quadrature algorithm is executed with the user-provided driving polynomial(s) (in full dimensional space) with masks identically equal to $1$ (as an easy means to convey that any subcell may contain the interface). In \cref{algo:main}, a mask $m$ is deemed ``nonempty'' if and only if $m(i) = 1$ on at least one subcell $i$.
\item Given a polynomial $f$ and mask $f_m$, the routine $\verb|nonzeroMask|(f,f_m)$
uses ``orthant tests'' to return a new mask $m$ such that $m(i) = 0$ if and only if at least one of the following conditions hold: (i) $f_m(i) = 0$; or (ii) it can be proven that $f$ does not have any real roots in the (slightly expanded) subcell $\Omega_i \pm \epsilon$. This routine can be used to determine which subcells contain the interface of $f$; in particular, if the resulting mask is empty, i.e., identically zero, then the polynomial $f$ can be removed from the dimension reduction process.
\item Given two polynomials $f$ and $g$, with masks $f_m$ and $g_m$, the routine $\verb|intersectionMask|\bigl((f,f_m),(g,g_m)\bigr)$
uses orthant tests to return a new mask $m$ such that $m(i) = 0$ if and only if at least one of the following three conditions hold: (i) $f_m(i) = 0$; or (ii) $g_m(i) = 0$; or (iii) it can be proven that $f$ and $g$ do not share any real roots in the (slightly expanded) subcell $\Omega_i \pm \epsilon$. This routine is used to quickly ascertain whether two polynomials have nonintersecting zero level sets as well as to help eliminate unwanted zeros of resultant/discriminant calculations, as explained in \cref{sec:resultants}. 

\item The routine $\verb|collapseMask|(m, k)$
takes a $d$-dimensional mask and collapses it to a $(d-1)$-dimensional mask by bitwise-\textsc{or}-ing along each column in the direction $k$, i.e., it returns $m_c : {\mathbb N}_{\mathcal M}^{d-1} \to \{0,1\}$ such that $m_c(i) = 0$ if and only if $m(i + \ell e_k) = 0$ for every $\ell = 1,\ldots,\mathcal M$. This routine, along with the related method \verb|faceRestriction| (below), is a kind of dimension-reduction analogue for masks: in \cref{algo:main}, these routines create a mask for the polynomials making up the domain of the base integral. 
\item Finally, the routine $\verb|faceRestriction|(m, k, y)$
takes a $d$-dimensional mask and returns its restriction to the lower ($y = 0$) or upper ($y = 1$) face in the direction $k$; in other words, the routine returns the mask $m_\text{face} : {\mathbb N}_{\mathcal M}^{d-1} \to \{0,1\}$ where $m_\text{face}(i) = m\bigl(i + [1 + y ({\mathcal M} - 1)] e_k \bigr)$.
\end{itemize}

The number of subcells used to compute these masks, $\mathcal M$, is a user-chosen parameter. With increasing $\mathcal M$, the accuracy of the associated Bernstein polynomial range evaluation increases, and thus so too the degree of confidence in isolating zeros of the involved polynomials. However, our experiments indicate that there is little benefit in employing a large $\mathcal M$, as a small $\mathcal M$ is often sufficiently effective in achieving the objectives of the masking operations. In particular, a value of 4 or 8 works well, and $\mathcal M = 8$ has been used throughout this work.\footnote{Profiling indicates that, for this choice of $\mathcal M$, mask computations represent a small to medium portion of the overall computational cost of the quadrature algorithm.} It is important to stress, however, that the ``correctness'' of the quadrature algorithms is independent of $\mathcal M$, meaning that the correct topology of the interface will always be determined by \cref{algo:main,algo:integrand}, even if no masks are used. The main benefit of the masks is to eliminate unwanted zeros of resultant/discriminant calculations, which results in less unnecessary interval subdivisions in \cref{algo:integrand}. Some of these aspects are discussed next.

\subsection{Resultants and pseudo-discriminants}
\label{sec:resultants}

Recall that the domain of the base integral is divided into pieces corresponding to the connected components of $(0,1)^{d-1} \setminus \bigl(\mathcal P \cup \mathcal Q \cup \mathcal  R\bigr)$, where (see \cref{fig:examples}): $\mathcal Q$ indicate the base points of upper- or lower-face crossings of the interface; $\mathcal R$ indicate the base points above which two interfaces cross; and $\mathcal P$ indicate base points above which the height function branches, typically at the vertical tangents of rounded surfaces. In all three cases, these sets can be implicitly defined by multivariate polynomials, with $\mathcal Q$ being a simple case of restricting the input polynomial(s) to the upper- and lower-faces. In the cases of $\mathcal P$ and $\mathcal R$, these sets can be implicitly defined through the use of polynomial resultant methods, as follows.

\subsubsection{Resultants of univariate Bernstein polynomials} 

Given two univariate polynomials $f$ and $g$, a \textit{resultant} $R(f,g)$ is a polynomial expression in the coefficients of $f$ and $g$ such that $R(f,g) = 0$ if and only if $f$ and $g$ share a common root (whether real or complex). A number of resultant formulations are possible; in many cases they are given by determinants of a Sylvester or B\'ezout matrix. These matrices are often stated in terms of the power/monomial basis coefficients of $f$ and $g$, however a number of recent works \cite{WINKLER2000179,WINKLER2003153,BINI2004319,WINKLER2008529} have derived formulations for use with the Bernstein basis:
\begin{itemize}
\item \textit{Sylvester resultant for Bernstein polynomials.} Here, we adopt the formulation given in \cite{WINKLER2008529}: the Sylvester resultant matrix, of size $(n + m) \times (n + m)$, corresponding to degree $n$ and $m$ polynomials $f(x) = \sum_{i=0}^n f_i b_i^n(x)$ and $g(x) = \sum_{i=0}^m g_i b_i^m(x)$, is given by
\[ S(f,g) = \begingroup %
\setlength\arraycolsep{4pt} \begin{bmatrix}
 f_0 \binom{n}{0} & f_1 \binom{n}{1} & f_2 \binom{n}{2} & \cdots \\
                  & f_0 \binom{n}{0} & f_1 \binom{n}{1} & f_2 \binom{n}{2} & \cdots\vphantom{\diagdots} \\
									&                  & f_0 \binom{n}{0} & f_1 \binom{n}{1} & \diagdots  \\
									& & & \diagdots & \diagdots \\
 g_0 \binom{m}{0} & g_1 \binom{m}{1} & g_2 \binom{m}{2} & \cdots\vphantom{\diagdots} \\
                  & g_0 \binom{m}{0} & g_1 \binom{m}{1} & g_2 \binom{m}{2} & \cdots\vphantom{\diagdots} \\
									&                  & g_0 \binom{m}{0} & g_1 \binom{m}{1} & \diagdots  \\
									& & & \diagdots & \diagdots 
\end{bmatrix} D^{-1}, \endgroup
 \]
where $D = \text{diag}\bigl( \{\binom{m+n-1}{i}\}_{i=0}^{m+n-1} \bigr)$ is a diagonal scaling.\footnote{For an intuitive derivation of the Sylvester matrix $S(f,g)$ for Bernstein polynomials, see \cite{FAROUKI19881}; this derivation does not include the diagonal scaling $D$, which can be effective at preventing an unnecessary rapid growth of the matrix entries for large degree polynomials, improving numerical conditioning.} Each row of $S(f,g)$ contains a shifted copy of the scaled Bernstein coefficients of $f$ and $g$, with the first $m$ rows corresponding to $f$ and the next $n$ rows corresponding to $g$. %
Using the Sylvester matrix, a resultant expression is then given by
\[ R(f,g) = \det S(f,g). \]

\item \textit{B\'ezout resultant matrix for Bernstein polynomials.} The (symmetric) B\'ezout resultant matrix $B(f,g) = ({\mathsf b}_{ij})_{i,j=1}^n$, of size $n \times n$, corresponding to degree $n$ polynomials $f(x) = \sum_{i=0}^n f_i b_i^n(x)$ and $g(x) = \sum_{i=0}^n g_i b_i^n(x)$ can be defined recursively\footnote{In floating point arithmetic, it is important to first build the lower left triangular portion of $B(f,g)$, and then copy the result (by symmetry) into the upper right triangular portion; a direct application of the recursive formula to compute the upper right triangular portion is numerically unstable and can result in a blowup of rounding errors.} \cite{BINI2004319} such that
\begin{align*}
 {\mathsf b}_{i,1} &= \smash{\tfrac{n}{i}} (f_i g_0 - f_0 g_i), &  1 \leq i \leq n , \\
 {\mathsf b}_{n,j+1} &= \smash{\tfrac{n}{n - j}} (f_n g_j - f_j g_n), &  \quad 1 \leq j \leq n - 1, \\
 {\mathsf b}_{i,j+1} &= \smash{\tfrac{n^2}{i(n - j)}} ( f_i g_j - f_j g_i ) + \smash{\tfrac{ j(n - i) }{i(n - j)}} {\mathsf b}_{i+1,j}, & 1 \leq i, j \leq n - 1.
\end{align*}
As before, a resultant expression is then given by
\[ R(f,g) = \det B(f,g). \]
\end{itemize}
Comparing these two approaches, the B\'ezout method has a number of potentially useful properties including, e.g., it results in a smaller matrix, and thus a faster determinant calculation; it can also exhibit better numerical conditioning than the Sylvester method \cite{WINKLER2013410}. However, it requires equal degree polynomials, thus requiring degree elevation of one of the polynomials in the mixed degree case, whereas the Sylvester method is directly applicable to mixed degrees. In this work, we have used the B\'ezout method whenever the input polynomials have equal degree, and the Sylvester method otherwise.

\subsubsection{Resultants of multivariate Bernstein polynomials} 
\label{subsec:resultants}

The above univariate resultant methods can be adapted to multivariate Bernstein polynomials, as follows. Let $f(x) = \smash{\sum_i f_i \prod_{\ell=1}^d b_{i_\ell}^{n_\ell} (x_\ell)}$ and $g(x) = \smash{\sum_i g_i \prod_{\ell=1}^d b_{i_\ell}^{m_\ell} (x_\ell)}$ denote two tensor-product Bernstein polynomials of degrees $(n_1,\ldots,n_d)$ and $(m_1,\ldots,m_d)$, respectively. Suppose that $k \in \{1, \ldots, d\}$ is a chosen elimination axis. Then, $f$ and $g$ can be reinterpreted as univariate polynomials of degree $n_k$ and $m_k$, with their coefficients a tensor-product Bernstein polynomial in one fewer dimensions. Using the lifting notation defined earlier, for a base point $\bar x = (x_1, \ldots, x_{k-1},x_{k+1},\ldots,x_d)$, denote these univariate polynomials by $f_{\bar x}$ and $g_{\bar x}$, such that $f_{\bar x}(y) = f(\bar x + y e_k)$ for all $y$, and similarly for $g_{\bar x}$. The above univariate resultant method can then be applied to $f_{\bar x}$ and $g_{\bar x}$ to obtain a resultant $R(f_{\bar x}, g_{\bar x})$. Importantly, the output of this operation is a polynomial in the coefficients of $f_{\bar x}$ and $g_{\bar x}$, and thus a polynomial in $\bar x$; denoting this polynomial by $R(f,g; k)$ with the argument $k$ indicating which variable or axis has been eliminated, we have
\begin{equation} \label{eq:resultant1d} \bigl[R(f, g; k)\bigr] (\bar x) = R(f_{\bar x}, g_{\bar x}) \text{ for all } \bar x. \end{equation}
For an arbitrary pair of tensor-product polynomials $f$ and $g$, the resultant $R(f,g;k)$ will also be a tensor-product polynomial. The maximum possible degree of this polynomial, $(r_1, \ldots, r_{k-1}, r_{k+1}, \ldots, r_d)$, is given by
\begin{equation} \label{eq:r} r_\ell = n_\ell m_k + m_\ell n_k, \quad \text{for} \quad \ell = 1, \ldots, k-1, k+1, \ldots, d. \end{equation}
A number of different approaches could be used to compute $R(f,g;k)$. For example, one could apply a division-free determinant algorithm together with Bernstein multiplication, to compute the determinant of the Sylvester or B\'ezout matrix for all $\bar x$ at once. In this work, we have instead opted for a faster method (in terms of complexity): the method uses a Bernstein interpolation algorithm applied to a pointwise nodal evaluation of $R(f_{\bar x}, g_{\bar x})$ at sufficiently many points to recover a tensor-product Bernstein polynomial of degree given by \eqref{eq:r}. Details on this algorithm are given in the supplementary material.

\subsubsection{Pseudo-discriminants}

Similar to the resultant, another key tool in algebraic geometry is the \textit{discriminant}. Given a polynomial $f$, the discriminant $\Delta(f)$ is a polynomial expression in the coefficients of $f$ such that $\Delta(f) = 0$ if and only if $f$ has a double root. One possible method for defining the discriminant is by a normalisation of the resultant of $f$ and its derivative, with the specific kind of normalisation sometimes varying from author to author. In this work, however, we deliberately choose to leave the expression unnormalised---in particular, we define the \textit{pseudo-discriminant} $\dotDelta$ of a multivariate polynomial $f$ in the direction $k$ as
\[ \dotDelta(f; k) := R(f, \partial_k f; k). \]
To illustrate the usefulness of this choice, consider the 2D polynomial $\phi(x,y) = xy$, whose zero level set is a cross aligning with the Cartesian coordinate axes. For the purposes of quadrature, supposing the height direction is $\hat y$, say, it is important to detect the location of the vertical interface, so that an accurate quadrature scheme can be built for all four quadrants. To do so, we can use the pseudo-discriminant of $\phi$ in the direction $\hat y$, whose calculation yields the 1D polynomial $\dotDelta(\phi; y) = R(\phi, \partial_y \phi; y) = x$; note that the zero set of this polynomial successfully pinpoints the vertical interface of $\phi$. On the other hand, the conventional discriminant of $\phi$ (using one of the most common kinds of normalisation) with elimination variable $y$ yields $\Delta(\phi; y) = 1$; as such, it does not provide any useful information about the position of the vertical interface. In a sense, a vertical interface is a degenerate example of a branching point (labelled ${\mathcal P}$ in \cref{fig:examples}); the conventional discriminant may not always pick up these points, but the pseudo-discriminant can, hence its adoption in the quadrature algorithms of this work. To compute $\dotDelta(f; k)$, the same algorithm as for the resultant can be used; for an arbitrary tensor-product polynomial $f$ of Bernstein degree $(n_1,\ldots,n_d)$, the output will also be a tensor-product polynomial, and its maximum possible degree $(r_1, \ldots, r_{k-1}, r_{k+1}, \ldots, r_d)$ is given by
\begin{equation} \label{eq:d} r_\ell = (2 n_k - 1) n_\ell, \quad \text{for} \quad \ell = 1, \ldots, k-1, k+1, \ldots, d. \end{equation}

\subsubsection{Application to the quadrature algorithm} \label{sec:spurious} Having established these tools to compute multivariate resultants and pseudo-discriminants, we are now in a position to describe their use in the quadrature algorithm. 

First, we consider the points $\mathcal R$ in \cref{fig:examples}, i.e., the base points above which the zero level sets of two polynomials $\phi_i$ and $\phi_j$ intersect. To capture the location of these points, we can make use of the resultant of $\phi_i$ and $\phi_j$ in the height direction $e_k$. In particular, given a base point $x \in (0,1)^{d-1}$, let $\phi_{i,x}$ denote the univariate polynomial corresponding to the vertical line restriction of $\phi_i$ at $x$, i.e., so that $\phi_{i,x}(y) = \phi_i(x + y e_k)$ for all $y$, and similarly for $\phi_j$. Then, $\bigl[R(\phi_i, \phi_j; k)\bigr](x) = 0$ if and only if $\phi_{i,x}$ and $\phi_{j,x}$ share a common root (real or complex). Consequently, the set of points ${\mathcal R}$ corresponding to the interfacial intersection of $\phi_i$ and $\phi_j$ is a subset of the zero isosurface of the polynomial $R(\phi_i, \phi_j; k)$.

A similar observation holds for the base points above which interfacial surfaces branch, i.e., where the multi-valued height function has a double root (points $\mathcal P$ in \cref{fig:examples}). We can capture the location of these points by making use of the pseudo-discriminant of $\phi_i$. In particular, $\bigl[\dotDelta(\phi_i;k)\bigr](x) = 0$ if and only if $\phi_{i,x}$ and its derivative share a common root (real or complex). %
As such, the set of ${\mathcal P}$ points corresponding to $\phi_i$ is a subset of the zero isosurface of the polynomial $\dotDelta(\phi_i; k)$.

\begin{figure}%
\centering\includegraphics{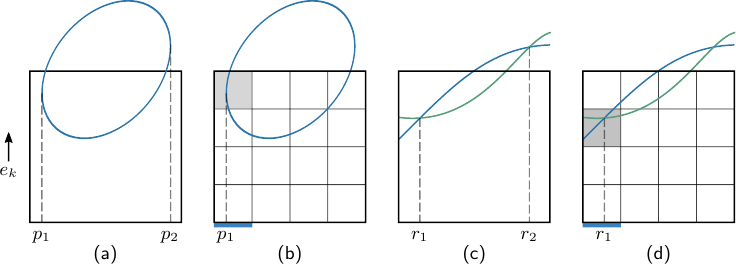}%
\caption{Superfluous zeros of resultant calculations and their mitigation through masking methods. \textit{(a)} The discriminant of a rotated ellipse yields a univariate polynomial with two real roots, $p_1$ and $p_2$; in particular, $p_2$ corresponds to a branching point of the interface located outside the domain of interest. (b) Using the masking methods of \cref{sec:masks}, here illustrated with $\mathcal M = 4$, the discriminant polynomial is associated with a mask which is enabled (value $1$) in the shaded cell, and disabled (value $0$) everywhere else. The restriction of this mask to the base domain yields a one-dimensional mask with one enabled cell (blue). Combined, the two masks keep $p_1$ and filters out the superfluous zero $p_2$. \textit{(c)} A resultant calculation of two polynomials, each implicitly defining the blue and green curves, yields a univariate polynomial with two real roots, $r_1$ and $r_2$; in particular, $r_2$ corresponds to a crossing point located outside the domain of interest. \textit{(d)} Similar to the left example, mask computations effectively filter out the superfluous root $r_2$.}
\label{fig:spurious}
\end{figure}

Note that, in both cases, the sets $\mathcal R$ and $\mathcal P$ may be proper subsets of the zero level sets of the respective resultant and pseudo-discriminant calculations. The reason for this is that resultant methods detect any shared roots, whether real or complex, including shared real roots that are outside the hyperrectangle of main interest to the quadrature algorithm. Some examples are shown in \cref{fig:spurious}. In \cref{fig:spurious}(a), the (pseudo-)discriminant calculation for a rotated ellipse yields a polynomial which has two real roots in the interval $\{0 \leq x \leq 1\}$, but only one of these roots corresponds to a branching point of the ellipse inside the unit square, whereas the other root is ``superfluous'' in the sense it detects a branching point outside the square of interest. In \cref{fig:spurious}(c), the resultant calculation of two intersecting curves correctly identifies the base point of one junction, but also produces a superfluous base point for a junction that occurs outside the square domain.

These superfluous zeros produced by resultant and pseudo-discriminant calculations do not affect the ``correctness'' of the quadrature algorithms, but they can cause unnecessary subdivision in subsequent one-dimensional quadrature calculations. To help mitigate these effects, and in many cases eliminate them altogether, we make use of the masking operations described in \cref{sec:masks}. Specifically, the routine \texttt{intersectionMask} is used to produce a mask whose subcells are marked ``$0$'' if it can be provably determined that the two input polynomials do not share any common roots. These masks are then associated with corresponding resultant or pseudo-discriminant polynomials and essentially state that any roots of these polynomials in a subcell marked ``$0$'' can be ignored. For the examples of \cref{fig:spurious}, these masks are illustrated in \cref{fig:spurious}(b,d), wherein the non-shaded subcells correspond to the disabled parts of the mask; a root of the 1D resultant or pseudo-discriminant in one of the non-shaded subcells is then ignored by the quadrature algorithm. These masking operations are also effective in 3D. Furthermore, they are often effective in proving the non-existence of any $\mathcal R$ or $\mathcal P$ points, in which case the associated resultant or pseudo-discriminant calculations can be avoided altogether.

In general, the above methods for computing resultants and pseudo-discriminants are sufficient for the vast majority of applications we have in mind for the quadrature algorithms developed in this work. In exceptional cases, or deliberately contrived situations, additional tools may be needed to fully realise these methods. For example, if the user-provided multivariate polynomial is the square of another polynomial, it may need to be made square-free by using GCD factoring methods.\footnote{Note, however, implicitly defined geometry given by the square of a polynomial does not lead to a well-posed problem for quadrature purposes: perturbations of the squared polynomial can result in the entire removal of the interface.} It is possible to incorporate additional steps into \cref{algo:main} to automatically handle these edge cases, however we have omitted these details here for the benefit of clarity. 

\subsection{Choosing a height direction}
\label{sec:height}

At each stage of the dimension reduction process, a height direction/elimination axis $e_k$ must be chosen. So far, we have seen how this choice affects the geometry of the base integral's domain via characterisation of the points $\mathcal P$, $\mathcal Q$, and $\mathcal R$ illustrated in \cref{fig:examples}, whose calculation, in turn, can depend on resultant and pseudo-discriminant calculations, as discussed in \cref{sec:resultants}. The choice of height direction also impacts the accuracy of the subsequent quadrature schemes. Some of these effects include:
\begin{enumerate}[label=(\roman*)]
\item Resultant and pseudo-discriminant computations can produce high degree polynomials in the description of the base integral's domain, an effect which can potentially compound further down the dimension reduction process for 3D (and higher dimension) problems. As such, it can be beneficial to limit the application of resultant and/or pseudo-discriminant polynomials.
\item The polynomials representing upper- and lower-face restrictions, i.e., those which characterise the points $\mathcal Q$, do not increase the degree of the input polynomial. As such, their creation and manipulation have a relatively minor impact on the overall complexity. However, if some faces of the hyperrectangle of $(0,1)^d$ have fewer or no interfacial intersections compared to other faces, it is preferential to use the corresponding height direction, as this means the geometry characterising the base integral's domain could be simpler, or even trivial.
\item As will be seen in the results, the accuracy of the overall quadrature scheme benefits from an effective application of Gauss-Legendre quadrature; in turn, these schemes work best if they are applied to interfacial geometry characterised by smooth, well-behaved height functions. In particular, a height direction in which there are no vertical tangents in the interface, i.e., $\mathcal P$ is empty, is preferred.
\end{enumerate}
These objectives are not always mutually compatible. For example, a direction which avoids vertical tangents may necessitate non-trivial upper or lower face restrictions. However, a clear trend found in this work is that the best quadrature results---in terms of both accuracy and efficiency of construction---are obtained if $e_k$ is chosen to be in the direction ``most'' orthogonal to the interfacial geometry. This metric prioritises the objectives (i) and (iii) above, in that it avoids vertical tangents and pseudo-discriminant calculations and greatly increases the odds of attaining high accuracy through effective application of Gauss-Legendre quadrature.

This heuristic is implemented follows. Suppose we are given a list of polynomials, $\{\phi_i\}$. First, we determine (cheaply, using masking operations) whether there are any directions $k$ such that the interface of every given polynomial is provably the graph of a (potentially multi-valued) height function devoid of vertical tangents/branches.\footnote{To do so, we can apply the \texttt{intersectionMask} routine to the polynomials $\phi_i$ and $\partial_k \phi_i$, and test whether the resulting mask is empty, i.e., identically $0$-valued.} Let $\mathcal K$ denote the set of all such directions $k$; if $\mathcal K$ is nonempty, then only those directions shall be considered in subsequent steps; otherwise, we set $\mathcal K = \{1, \ldots, d\}$ and consider all directions. 

The height direction in \cref{algo:main} on line \ref{algo:mainchoice} is then defined by choosing $k \in \mathcal K$ with the best ``score'', as follows (if $\mathcal K$ is a singleton, this step is trivial). For each $i$, let $m_i \in {\mathbb N}_{\mathcal M}^d \to \{0,1\}$ be the mask indicating which subcells of $\phi_i$ are provably nonzero (as computed by, e.g., \texttt{nonzeroMask}, see \cref{sec:masks}). Then the score in direction $k$ is evaluated by
\[ s_k = \sum_i \sum_{\substack{\ell \in {\mathbb N}_{\mathcal M}^d \\ m_i(\ell) \neq 0}} \frac{|\partial_k \phi_i (x_\ell)|}{\|\nabla \phi_i (x_\ell) \|_1} \]
where $x_\ell$ is the midpoint of subcell $\ell$. In essence, this metric samples the gradient of $\phi_i$ at the midpoint of every mask subcell which (likely) contains the interface, weighted so that polynomials $\phi_i$ having greater interfacial measure have a greater influence in the overall score. In this fashion, the direction $k \in \mathcal K$ having largest score $s_k$ heuristically determines which coordinate direction is most orthogonal to the combined interface(s) of every $\phi_i$. Importantly, these mask operations and score evaluations are cheap and have negligible cost in the overall quadrature algorithm.

\subsection{One-dimensional quadrature}
\label{sec:1dquad}

One of the last remaining design aspects of the quadrature algorithm concerns the one-dimensional vertical line integrals, i.e., the choice of quadrature schemes in the implementation of \cref{algo:integrand}. Recall these line integrals correspond to evaluating, for a given base point $x \in (0,1)^{d-1}$, the base integral's integrand ${\mathcal F}(x)$, where ${\mathcal F}(x) = \circint_{(0,1) \setminus \Gamma(x)} f(x + y e_k)\,dy$, with $f$ being the inner integral's integrand, while $\Gamma(x) = \bigcup_{i=1}^n \{y : 0 \leq y \leq 1 \text{ and } \phi_i(x + y e_k) = 0 \}$ specifies the set of interface points above $x$.

To approximate ${\mathcal F}(x)$ for a given $x$, we first use root finding algorithms to compute the set of interfacial points $\Gamma(x)$. (Fast and backward stable methods for doing so are discussed in the supplementary material.) These roots are then used to partition $[0,1]$ into a set of one or more subintervals, after which a one-dimensional quadrature scheme is applied to each subinterval. With the goal of obtaining a simple ``black box'' quadrature algorithm, we here make two highly-simplifying design choices:
\begin{enumerate}[label=(\roman*)]
\item Every subinterval shall apply the same kind of quadrature scheme, with the same number of quadrature points $q$ per subinterval; and 
\item The choice of quadrature scheme, whether it be Gauss-Legendre or tanh-sinh or an alternative, shall be independent of the base point $x$; in other words, the choice is made on a per-dimension basis. %
\end{enumerate}
Clearly, this approach is not at all adaptive; e.g., choosing $q$ fixed can lead to an unnecessary number of quadrature points for small subintervals. The approach is also somewhat blind to the interfacial geometry at hand: different parts of $\Gamma$ may exhibit different kinds of geometry and thus benefit from different kinds of quadrature that could vary throughout the domain of the base integral; as such, the universal scheme required by (ii) above may not always yield the smallest errors. Our experiments show, however, that the above simplifications nevertheless result in highly accurate quadrature schemes which exhibit exponential convergence in $q$. Additional remarks about these choices, including the possibility of adaptive methods, are given in the concluding remarks.

In our experiments, we have mainly explored the use of Gauss-Legendre and tanh-sinh quadrature schemes. The application of Gauss-Legendre methods (see, e.g., \cite{doi:10.1137/120889873}) should be familiar to most readers and is not discussed further here; the use of tanh-sinh, however, is worthy of additional motivation. This scheme, also known as the double-exponential method, was originally developed by Takahasi and Mori \cite{takahasi1974double} and makes use of the transformation $x = \tanh( \tfrac12 \pi \sinh(\xi))$ to map the interval $x \in (-1,1)$ to the whole real line $\xi \in (-\infty,\infty)$ where then a trapezoidal rule is applied. The method can be highly effective at handling endpoint singularities, where in many cases it can achieve (near) exponential convergence. In the present work, our motivation for its application is to treat the cases of vertical tangents in the multi-valued height function; these vertical tangents give rise to an endpoint singularity in the base integral's integrand ${\mathcal F}$. As an example, consider the area integral of $f$ on the upper half circle with height function $h(x) = \smash{\sqrt{1 - x^2}}$; the task of the base integral is to then compute $\smash{\int_{-1}^1 \bigl( \int_0^{h(x)} f(x, y)\,dy \bigr)\,dx}$; for $f \equiv 1$, in particular, this reduces to $\smash{\int_{-1}^1 \sqrt{1 - x^2}\,dx}$. Applying Gauss-Legendre to this integral yields a poor convergence rate as $q$ increases; on the other hand, tanh-sinh yields rapid, near exponential convergence \cite{doi:10.1137/130932132}.

Relative to its standard reference domain, $(-1,1)$, the tanh-sinh quadrature method yields a quadrature scheme of the form $\smash{\int_{-1}^1 f(x)\,dx \approx \sum_\ell w_\ell f(x_\ell)}$; the standard tanh-sinh schemes usually have an odd number of quadrature points and do not exactly integrate constant functions, i.e., $\sum_\ell w_\ell \neq 2$. In the present work, we have implemented small modifications to allow for an even number of quadrature points and normalise the weights so that they exactly integrate constant functions;\footnote{This normalisation results in the guarantee that any quadrature scheme produced by the dimension reduction algorithm satisfies the condition $\int_{(0,1)^d \setminus \Gamma} 1 = 1$.} importantly, these modifications do not alter the exponential convergence of these methods. Details are given in the supplementary material.

To summarise, \cref{algo:integrand} approximates the vertical line integrals ${\mathcal F}(x)$ in the dimension reduction quadrature algorithm. Each subinterval uses, for simplicity, the same number of quadrature points $q$; each dimension shall apply, for simplicity, a fixed choice of Gauss-Legendre or tanh-sinh quadrature schemes, independent of the argument $x$. As such, the user may choose $q$ and may choose Gauss-Legendre or tanh-sinh for each dimension of the reduction process. These choices, particularly of whether to apply Gauss-Legendre or tanh-sinh as a function of the detected interfacial geometry, is one of the main topics of discussion for the experimental results presented in \cref{sec:results}. %

\subsection{Surface integration}
\label{sec:surface}

In the formulation so far, we have only considered the case of volumetric integrals over implicitly defined domains. In essence, the same dimension reduction algorithm can be used in the case of surface integrals---the main modification needed is to replace the inner-most vertical line integral (corresponding to the very first elimination axis) with a new functional which evaluates the user-provided integrand $f$ at the surface, appropriately weighted by Jacobian factors to account for its curvature. However, in some cases of interfacial geometry, particularly degenerate cases, it is important to appropriately handle singular Jacobian factors; here, we discuss an approach which robustly handles these situations.

The problem set up is as follows. Suppose we are given one or more polynomial functions $\phi_i : (0,1)^d \to \R$, $i = 1, \ldots, n$, which together implicitly define an interface $\Gamma := \bigcup_{i=1}^n \{\phi_i = 0\}$. Suppose also that we are given an integrand $f : \Gamma \to \R$, which need only be defined on the interface $\Gamma$. Our task is to build a high order quadrature scheme approximating $\int_\Gamma f$. Applying the same dimension reduction ideas as in \cref{sec:mainformulation}, suppose we have chosen a height direction, $e_k$, say. Then, we may reinterpret $\Gamma$ as the graph of a (potentially multi-valued) height function and compute the surface integral correspondingly. The curvature of the height function induces a surface Jacobian factor: in \cite{HighorderImplicitQuad} it is shown (via the implicit function theorem) this can be computed via the gradient of the associated level set function. Using these ideas, we obtain $\int_\Gamma f = \int_\Omega {\mathcal F}(x)\, dx$ where $\Omega$ is the domain of the base integral, $x \in \Omega \subseteq (0,1)^{d-1}$ is a base point, and
\begin{equation} \label{eq:surf1} {\mathcal F}(x) = \sum_{i = 1}^n \sum_{y \in \Gamma_i(x)} f(x + y e_k) \left. \frac{ |\nabla \phi_i| }{ |\partial_k \phi_i| } \right|_{x + y e_k} \end{equation}
where
\[  \Gamma_i(x) := \{y : 0 \leq y \leq 1 \text{ and } \phi_i(x + y e_k) = 0 \}. \]
In essence, the base integral's integrand ${\mathcal F}(x)$ evaluates the given integrand $f$ at every interfacial point above $x$; at each such point, the curvature of the interface is accounted for by the Jacobian factor $|\nabla \phi|/|\partial_k \phi|$.

If there is a height direction $e_k$ such that every such Jacobian factor is well-defined, we can then employ the existing dimension reduction algorithm---in doing so, a quadrature scheme for $\int_\Gamma f = \circint_\Omega {\mathcal F}(x)\,dx$ is constructed, with $\Omega$ defined in the same way as the volumetric case, i.e., as the set of points $(0,1)^{d-1} \setminus (\mathcal Q \cup \mathcal R)$ in which (see \cref{fig:examples}) the upper- and lower-face interfacial crossings $\mathcal Q$ and internal junctions $\mathcal R$ have been removed (there are no points $\mathcal P$ by assumption). However, this approach cannot be applied if no such height direction exists.

Two examples help to illustrate the situation at hand, the second more severe than the first:
\begin{enumerate}[label=(\roman*)]
\item Consider $\phi$ given by $\phi(x,y) = (x - \tfrac12)^2 + (y - \tfrac12)^2 - \tfrac{1}{16}$, whose interface is a circle of radius $0.25$ centred at $(0.5,0.5)$. Then, every elimination axis yields a height function with vertical tangents; e.g., for $e_k = \hat y$, the height function %
has a vertical tangent at $x \in \{0.25,0.75\}$. At these locations, the surface Jacobian factor $|\nabla \phi|/|\partial_y \phi| $ is not well-defined. One potential workaround is to apply tanh-sinh quadrature in the base integral and potentially yield exponential convergence;\footnote{One can show that, on the upper half of the circle with corresponding height function $h(x)$, the base integral's integrand takes the form ${\mathcal F}(x) = f(x,h(x)) \sqrt{1 + (h'(x))^2}$, where the latter term has a singularity of the form $1/\sqrt{\tilde x}$ with $\tilde x$ being the distance to the circle endpoints. A square-root singularity like this can be easily handled by tanh-sinh quadrature, provided certain measures are taken to guard against floating point roundoff errors \cite{baileynotes,vanherck2020tanhsinh}.} however, in more general cases it can be problematic to yield full machine precision due to the way tanh-sinh sends quadrature points exponentially close to the integral endpoints, combined with roundoff errors in computing polynomial roots.

\item  Consider $\phi$ given by $\phi(x,y) = (x - \tfrac12)(y - \tfrac12)$, whose interface is a cross of the form $\smash{\scalebox{1.5}{+}}$ centred at $(0.5,0.5)$. The approach used by \eqref{eq:surf1} is not all amenable to this example. For example, with elimination axis $\hat y$, the dimension reduction algorithm will correctly sense the location of the vertical arm of the interface, but the associated integrand ${\mathcal F}(x)$ given by \eqref{eq:surf1} at $x = 0$ is not at all well-defined. There is no means for the base integral to suitably place quadrature points on the vertical arm of the interface and to assign well-posed quadrature weights. However, the base integral can correctly handle the horizontal arm---this gives a hint on how to handle the general case, which is to apply the dimension reduction algorithm in all dimensions $e_k$, $k = 1, \ldots, d$, and then aggregate the results.
\end{enumerate}
Exploring this idea further, note that the surface Jacobian factor $|\nabla \phi|/|\partial_y \phi| $ in \eqref{eq:surf1} is precisely $1/|n_k|$, where $\mathbf n$ is the unit-magnitude normal to the interface $\Gamma$. As such, performing a surface integral, not of $f$, but of an integrand of the form $f n_k$, will cancel the problematic term. For such an integrand, ${\mathcal F}$ now takes the form
\begin{equation} \label{eq:surf2} {\mathcal F}(x) = \sum_{i = 1}^n \sum_{y \in \Gamma_i(x)} f(x + y e_k) \left.  \text{sign} \bigl( \partial_k \phi_i \bigr) \right|_{x + y e_k}. \end{equation}
This integrand is piecewise smooth on the connected components of the base integral's domain $\Omega = (0,1)^{d-1} \setminus \{\mathcal P \cup \mathcal Q \cup \mathcal R\}$; in particular the sign term is constant on each connected component due to the removal of the set $\mathcal P$ (recall that $\mathcal P$ characterises the base points above which the zero isosurfaces of $\phi_i$ and $\partial_k \phi_i$ intersect). As such, we can use the existing dimension reduction algorithm to build a high order quadrature scheme for $\int_\Gamma f n_k = \circint_\Omega {\mathcal F}(x)\,dx$. Performing this process for every height direction $k = 1, \ldots, d$ then yields a robust and accurate means to compute a quadrature scheme of $\int_\Gamma f \mathbf n$, for all sufficiently smooth $f$. 

In summary, the above considerations suggest the following two mechanisms: %
\begin{itemize}
\item If it can be proven that the interface is the graph of a (multi-valued) height function in direction $e_k$, without vertical tangents, and such that the associated Jacobian factors $|\nabla \phi|/|\partial_k \phi|$ are sufficiently well-behaved, then we can apply \cref{algo:main} where line \ref{algo:mainintegranddef} is replaced with a function evaluating \eqref{eq:surf1}.
\item Otherwise, apply an ``aggregated'' method. For each $k = 1, \ldots, d$, execute \cref{algo:main} with the following modifications: (i) overrule line \ref{algo:mainchoice} to use the height direction $k$; and (ii) replace line \ref{algo:mainintegranddef} with a suitable functional evaluating \eqref{eq:surf2}. For each $k$, the output of this process is a quadrature scheme approximating $\int_\Gamma f\,n_k$. 
\end{itemize}
Clearly, the aggregated method is computationally more expensive than the first approach, having to iterate over each possible choice of elimination axis. However it is more well-posed and less sensitive to perturbations in the input geometry.

Note that the aggregated method essentially produces a surface quadrature scheme in \textit{flux form}, i.e., given a sufficiently smooth function $f$, the algorithm produces a quadrature scheme $\int_\Gamma f \mathbf n \approx \sum_\ell {\mathbf w}_\ell f(x_\ell)$ where ${\mathbf w}_\ell$ are vector-valued quadrature weights. If needed, one could then build a surface quadrature scheme in \textit{non-flux form} by baking in the normal vector to yield a scheme of the form $\int_\Gamma g \approx \sum_\ell ({\mathbf w}_\ell \cdot {\mathbf n}_\ell) g(x_\ell)$, where ${\mathbf n}_\ell$ is the normal vector of the interface associated with $x_\ell$. %
Note that in many applications of surface quadrature, the flux form of the surface integral is the most relevant or natural. For example, in finite volume methods, discontinuous Galerkin methods, etc., the associated weak formulations require elemental boundary integrals which are typically in flux form, e.g., to compute $\int_\Gamma \hat f \mathbf n$ where $\hat f$ is a numerical flux. Our numerical experiments (see \cref{sec:results}) show that the flux-form of surface integration generally yields much higher accuracy than the non-flux form, essentially because the normal component $n_k$ exactly cancels the surface Jacobian factors $1/|n_k|$.\footnote{For example, in two dimensions, the flux-form integral of $f \equiv 1$, over the one-dimensional curve $\Gamma$, is computed exactly, independent of $q$.}

\subsection{Summary} %

In essence, the quadrature algorithm operates by recasting the implicitly defined interface as the graph of a multi-valued height function. Using the corresponding splitting of coordinates, the overall integral is reformulated as a base integral with an associated integrand ${\mathcal F}$. In the case of volumetric integrals, the role of ${\mathcal F}$ is to compute a vertical line integral in the height direction; in the case of surface integrals, the role of ${\mathcal F}$ is to evaluate the user-provided integrand on the interface $\Gamma$, appropriately modified by curved surface Jacobian factors. Continuing recursively, we are led to a conceptually simple dimension reduction algorithm---in order to achieve a high level of accuracy, it is important to appropriately identify and treat the piecewise smooth behaviour of ${\mathcal F}$. In turn, this corresponds to identifying the different characteristics of the height function, notably where $\Gamma$ crosses the upper- and lower-faces, its branching points, and the intersection points. The latter two require some workhorse tools from algebraic geometry, i.e., resultants and (pseudo-)discriminants, and are perhaps the most complex components of the quadrature algorithm. In \cref{sec:spurious}, we saw that these resultant methods can sometimes produce superfluous zeros, some of which correspond to branching/junction points outside the unit cube of interest. Many of these superfluous zeros can be removed by the masking methods. The masks also provide additional means for (cheaply) characterising the interface geometry, and we use them to help decide which coordinate axis to use as a height direction, as well as guide decisions on which kind of one-dimensional quadrature to use.

These ideas led to the main driver algorithms, \cref{algo:main} and \cref{algo:integrand}. Starting with the top-level integral problem, a coordinate axis is chosen for the height direction, and then a set of polynomials is constructed to determine the implicitly defined domain of the base integral. This process continues recursively, eliminating one coordinate axis at a time, down to the one-dimensional base case. The base case then performs a one-dimensional quadrature and this process continues recursively back up the tree: each integrand evaluation adds one more dimension at a time, terminating with the evaluation of the user-supplied top-level integrand $f$. %
In particular:

\begin{itemize}
\item To compute a volumetric quadrature scheme corresponding to
\[ \circint_{(0,1)^d \setminus \Gamma} f \approx \sum_i w_i f(x_i), \]
where $\Gamma = \bigcup_{i = 1}^n \{ \phi_i = 0 \}$ is implicitly defined by $n \geq 1$ user-defined polynomials $\phi_i$, \cref{algo:main} is executed with the integrand $f$ and every mask initialised by $m_i \equiv 1$. These mask values can be viewed as a simple default, conveying that any part of $(0,1)^d$ may contain the interface. The output is a quadrature scheme whose weights are, by construction, strictly positive. In fact, the nodes fill the unit cube and by normalisation conditions, we have $\sum_i w_i = 1$. Importantly, the scheme is accurate on each connected component of $(0,1)^d \setminus \Gamma$. As such, the sign of $\phi_i$ can be used to cluster the quadrature nodes; for example, to compute 
\[ \int_{(0,1)^d \cap \bigcap_{i=1}^n \{\phi_i < 0\}} f,\]
one can keep just the quadrature points satisfying $\phi(x_i) < 0$, leaving their weights unmodified. We note the algorithms can be easily adjusted if it was known ahead of time that only one side of $\Gamma$ is needed.

\item To compute a surface quadrature scheme, either for $\circint_{(0,1)^d \cap \Gamma} g$ or for $\circint_{(0,1)^d \cap \Gamma} g\mathbf n$, follow the advice of \cref{sec:surface}. The first is a surface integral in non-flux form, whereas the second is in flux-form and typically has higher accuracy. By construction, the resulting quadrature nodes are strictly on $(0,1)^d \cap \Gamma$. A similar operation as above could be used to cluster the nodes according to the different parts of a multi-polynomial interface.
\end{itemize}

\begin{figure}%
\centering\includegraphics{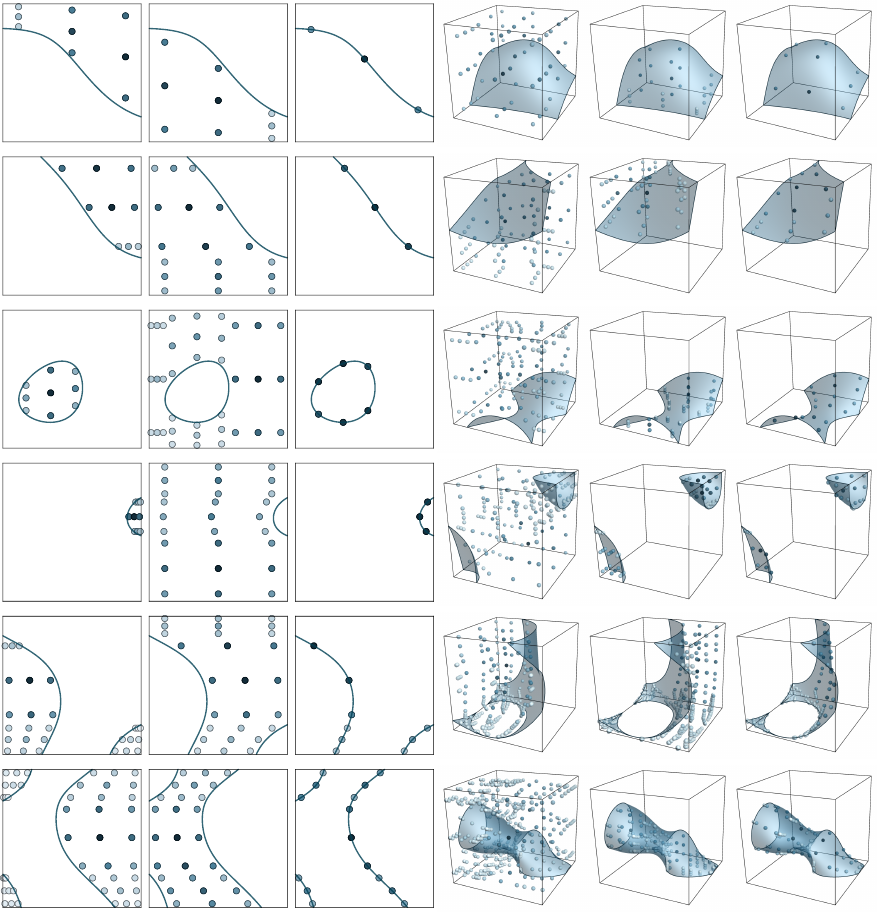}%
\caption{Examples illustrating the output of the quadrature algorithm for different kinds of interfacial geometry, using $q = 3$ quadrature nodes per one-dimensional integral. In each group of three, the left and middle figures correspond to volumetric integrals over the implicitly defined regions $\{\phi < 0\}$ and $\{\phi > 0\}$, while the right figure corresponds to the surface integral over $\{\phi = 0\}$; in each instance, the quadrature points are shaded according to their weight, such that low-valued weights are pale and higher-valued weights are darker.}%
\label{fig:quadexamples}%
\end{figure}

\cref{fig:quadexamples} illustrates the output of the quadrature algorithm for various kinds of interfacial geometry in 2D and 3D, with a single input polynomial (multiple polynomials are demonstrated later in \cref{sec:intersections}). In each case, the quadrature scheme consists of one or more patches of points in tensor-product form, with each patch warping to the geometry of the interface; this is most easily observed in the 2D examples and some of the simpler 3D examples. Naturally, the more complex 3D examples require more patches of points. 

\section{Numerical Experiments}
\label{sec:results}

In this section, we demonstrate the quadrature algorithm and investigate its accuracy through a variety of numerical experiments. The test problems range from interfacial geometry of mild to modest difficulty, yielding high-order convergence rates, through to more challenging scenarios involving cusps, pushing the limits of what can be achieved in finite precision arithmetic. In general, the quadrature schemes yield high degrees of accuracy, e.g., approximately exponential convergence in the number of quadrature points. In some cases, convergence to machine precision is quickly attained and this can prevent a clear observation of asymptotic convergence rates. %
To help with these effects, in some of our presented results we have made use of high-precision arithmetic via the QD library \cite{qdlib}, which provides double-double and quadruple-double floating point routines. We emphasize that high-precision arithmetic is not a requirement of the quadrature algorithms---these tools are only being used to help reveal the asymptotic convergence behaviour and quantify errors in relation to machine epsilon.

\subsection{Convergence under \texorpdfstring{$h$}{h}-refinement}
\label{sec:ellipse}

The majority of the presented results consider a fixed constraint cell $U$ and examine convergence under the action of $q$-refinement, i.e., as more and more quadrature points are used. In this study, however, we consider convergence under $h$-refinement, wherein both $q$ and the interfacial geometry is globally fixed and the cell $U$ is made smaller and smaller. Under this kind of refinement, for a globally smooth interface, the local interfacial geometry, $U \cap \Gamma$, becomes asymptotically flat; as such, we expect high rates of convergence, matching those of a Gaussian quadrature scheme. To confirm this, we consider a test problem of computing the surface area and volume of a $d$-dimensional ellipsoid:
\begin{itemize}
\item In 2D, the ellipse is given by the zero level set of $\phi(x,y) = x^2 + 4 y^2 - 1$; the ellipse has semimajor axis $1$, semiminor axis $\tfrac12$, an area of $V = \tfrac{\pi}{2}$, and perimeter $A = 4.84422411\ldots$.
\item In 3D, the ellipsoid is given by the zero level set of $\phi(x,y,z) = x^2 + 4 y^2 + 9z^2 - 1$; the ellipsoid has semiprincipal axes $1$, $\tfrac12$, and $\tfrac13$, a volume of $V = \tfrac{2\pi}{9}$, and surface area $A = 4.40080956\ldots$.
\end{itemize}
We consider a reference domain $(-1.1,1.1)^d$ enclosing the geometry, discretised by a uniform Cartesian grid consisting of $n$ cells in each direction, such that $h = 2.2/n$. For a fixed $q$, the quadrature algorithm is then applied to each cell, and the results summed over the entire grid so as to compute the total surface area and volume, denoted $A_{q,h}$ and $V_{q,h}$, respectively. In particular, and for simplicity of discussion, in these results we use Gauss-Legendre quadrature for all of the one-dimensional integrals executed by \cref{algo:integrand}.%

\begin{figure}%
\centering\includegraphics{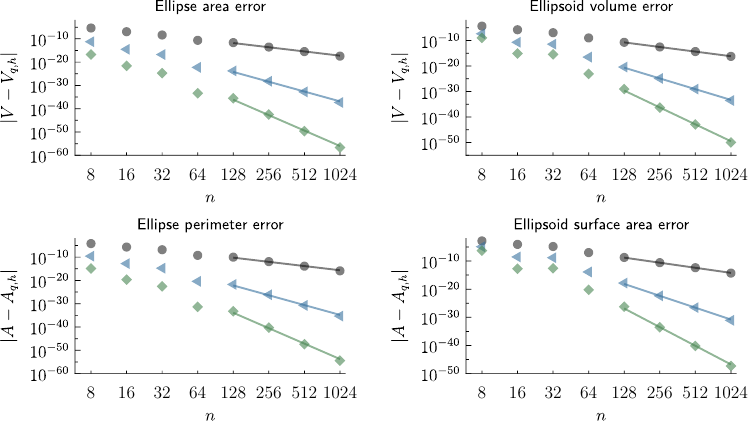}%
\caption{Results of the ellipsoidal test problems considered in \cref{sec:ellipse}, studying quadrature accuracy under the action of $h$-refinement, for three different choices of quadrature order $q$: $\bullet$ (black) denotes $q = 3$, $\trianglesymbol$ (blue) denotes $q = 7$, and $\squaresymbol$ (green) denotes $q = 11$; the corresponding lines have slope $2q$, indicating the ${\mathcal O}(h^{2q})$ asymptotic convergence rate as $h = 2.2/n \to 0$.}
\label{fig:ellipsehrefinement}
\end{figure}

In the presented results, we have made use of quadruple-double precision so as to more clearly reveal the asymptotic high-order convergence rates. \cref{fig:ellipsehrefinement} plots the errors in the computed surface area and volume of the ellipse and ellipsoid, for three different values of $q \in \{3,7,11\}$. In each case, the fitted line of slope $2q$ confirms that the quadrature schemes inherit the high-order convergence rate of Gaussian quadrature. For example, $22^\text{nd}$-order convergence rates are obtained with $q = 11$, on all four test problems. Similar results are obtained for other values of $q$: we have tested $q = 1, 2, \ldots, 11$ and observed ${\mathcal O}(h^{2q})$ errors in every case. Besides the ellipsoidal geometry of this example, the same conclusions hold as well for other kinds of sufficiently smooth geometry and integrand functionals. %

\subsection{Ellipsoid under \texorpdfstring{$q$}{q}-refinement}
\label{sec:ellipseQ}

\begin{figure}%
\centering\includegraphics{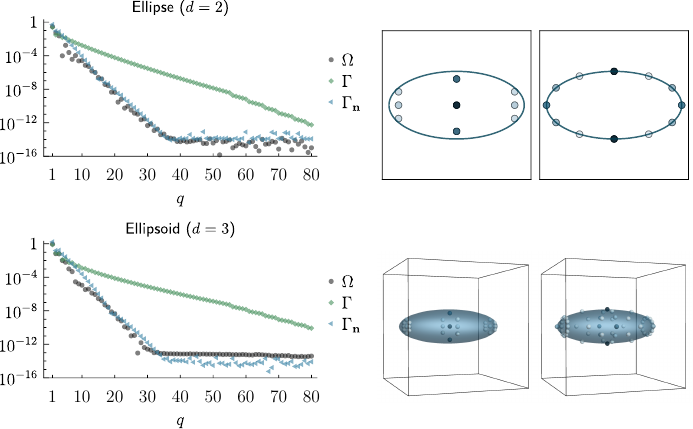}%
\caption{\textit{(left)} Quadrature accuracy as a function of $q$ for the ellipsoidal test problems considered in \cref{sec:ellipseQ}. Here, $\bullet$ (black) denotes the relative error in the computed volumetric integral, $I_\Omega$; $\squaresymbol$ (green) denotes the relative error in the computed surface integral in non-flux form, $I_\Gamma$; and $\trianglesymbol$ (blue) denotes the relative error in the computed flux-form surface integral, $I_{\Gamma_{\textbf n}}$. \textit{(right pair)} Illustrations of the $q = 3$ volumetric and (aggregated) surface quadrature schemes; in each instance, the quadrature points are shaded according to their weight, such that low-valued weights are pale and higher-valued weights are darker.}
\label{fig:ellipseqrefinement}
\end{figure}

Here we consider the same ellipsoidal geometry as in the previous problem, and examine quadrature accuracy under the action of $q$-refinement, i.e., keeping the geometry $U \cap \Gamma$ fixed while the quadrature order $q$ varies. In particular, we examine the error in computing the following three integral quantities,
\[ I_{\Omega} := \int_{U \cap \{\phi < 0\}} f, \quad I_{\Gamma} := \int_{U \cap \{\phi = 0\}} f, \quad I_{\Gamma_{\textbf n}} := \int_{U \cap \{\phi = 0\}} f \mathbf n, \]
where $U := (-1.1,1.1)^d$ is fixed, and the integrand $f : \R^d \to \R$ is given by\footnote{There is nothing particularly special about the chosen integrand $f$, other than being smooth, non-trivial, positive on $U$, and such that $I_{\Gamma_{\mathbf n}}$ is nonzero.} $f(x) = \cos (\tfrac14 \|x - \tfrac14\|_2^2)$. Here, $\phi$ denotes the level set function defining the ellipse/ellipsoid (identical that defined in \cref{sec:ellipse}), $I_\Omega$ denotes the volumetric integral of $f$ over the inside of the ellipse/ellipsoid, while $I_{\Gamma}$ and $I_{\Gamma_{\textbf n}}$ denote the surface integral of $f$ in non-flux and flux-form, respectively. Using the strategy outlined in \cref{sec:height} for choosing the height direction, in computing the  volumetric integral $I_\Omega$ the quadrature algorithm operates as follows:
\begin{itemize}
\item In 2D, the first elimination axis is $\hat y$, and so the outer integral is over the $x$-axis and the inner vertical line integrals occur in the $\hat y$ direction. The outer integral sees a height function with branching points (similar to \cref{fig:examples}, left); consequently, we apply tanh-sinh quadrature on the outer integral, and Gauss-Legendre quadrature on the inner integral. 
\item In 3D, the first elimination axis is $\hat z$, followed by elimination of the $\hat y$ axis, and so the outer-most integral is over $x$ and the inner-most integral over $z$. Similar to \cref{fig:examples}, left, the outer integrals involve height functions with vertical tangents, and so we apply tanh-sinh quadrature in \cref{algo:integrand} for integrals over the $x$ and $y$ direction, and apply Gauss-Legendre quadrature on the inner-most integral over $z$.
\end{itemize}
In both the 2D and 3D cases, the surface integrals are computed using the aggregated method discussed in \cref{sec:surface}. Illustrations of these schemes for $q = 3$ are shown in \cref{fig:ellipseqrefinement}.

For each integral, $I_{\Omega}$, $I_{\Gamma}$, and $I_{\Gamma_{\textbf n}}$, a reference solution (accurate to full machine precision) is computed using a high degree quadrature order, specifically $q = 100$. We then examine the error in these integrals, compared to the reference solution, using quadrature orders ranging from $q = 1$ up to $q = 80$; \cref{fig:ellipseqrefinement} contains the results,\footnote{This figure, and all others like it in this paper, plot the relative error compared to the reference/exact solution, e.g., $|I_{\Omega,q} - I_{\Omega,\text{exact}}| / |I_{\Omega,\text{exact}}|$. In the case of the flux-form surface integral, the vector maximum norm is used, i.e., $\|I_{\Gamma_{\textbf n},q} - I_{\Gamma_{\textbf n},\text{exact}}\|_\infty / \|I_{\Gamma_{\textbf n},\text{exact}}\|_\infty$.} computed using native double precision. In the case of the volumetric integral and the flux-form surface integral, we observe rapid,  approximately exponential convergence in $q$, with both integrals attaining maximum possible precision at around $q \approx 34 \pm 2$. This number is typical for tanh-sinh quadrature applied to integrals of the form $\int_{-1}^1 \sqrt{1 - x^2}$. However, for the surface integral in non-flux form, $I_{\Gamma}$, we observe a slower convergence rate, albeit also exponential. This slower convergence in computing $I_{\Gamma}$, versus $I_{\Gamma_{\textbf n}}$, confirms the discussion of \cref{sec:surface}; in particular, the method to compute $I_{\Gamma}$ is essentially computing the surface integral of $f |n_k|^2$, with the result summed over all $k$, where $\mathbf n$ is the normal to the interface. The corresponding integrand is bounded and smooth, but its analytic extension to the complex plane has a pole close enough to the origin that it impacts the convergence rate of tanh-sinh quadrature. For example, in the 2D case, one can show that, for elimination axis $\hat x$, the $y$ integral over $-\tfrac12 < y < \tfrac12$ has an integrand of the form $f(h(y),y) \sqrt{1 - 4 y^2} / \sqrt{1 + 12 y^2}$; this latter integrand has a pole inside the unit disc of the complex plane, and it is well-known that convergence of tanh-sinh slows in this case \cite{doi:10.1137/130932132}. In essence, the (aggregated) surface integral in non-flux form computes integrals with an extra factor of $|n_k|$ in the integrand (versus the flux-form), thereby changing its regularity, which in turn impacts the convergence rate. This phenomenon is observed in many of our results, and further comments are given in the concluding remarks. 

\subsection{Randomly generated geometry}
\label{sec:random}

We next consider quadrature performance on a class of randomly generated polynomials whose corresponding implicitly defined geometry exhibits a wide variety of characteristics, such as multi-component domains, multi-valued height functions, high degrees of curvature, and tunnels. Through the computation of aggregated statistics, our goal is to: (i) help determine an automated strategy for choosing between Gauss-Legendre and tanh-sinh quadrature in the overall quadrature algorithm; and (ii) assess quadrature convergence on the kinds of interfacial geometry which may arise in a complex computational fluid dynamics problem, e.g., modelling atomizing liquid.

The class of randomly generated polynomials is defined through the orthonormal Legendre polynomials, as follows. Let $p_i$, $i = 0, 1, 2$, denote the first three normalised univariate Legendre polynomials relative to the interval $(-1,1)$, i.e.,
\[ p_0(x) = \sqrt{\tfrac12}, \quad p_1(x) = \sqrt{\tfrac32}\, x, \quad p_2(x) = \sqrt{\tfrac58}  (3x^2 - 1). \]
In $d$ dimensions, define a tensor-product degree $(2,\ldots,2)$ polynomial $\phi: \R^d \to \R$ as follows:
\begin{equation} \label{eq:phipoly} \phi(x) = \sum_{i \in \{0,1,2\}^d} \lambda(i)\, c_i \prod_{\ell=1}^d p_{i_\ell} (x_\ell), \end{equation}
where $c_i \in \R$ is a given set of $3^d$ coefficients, one for each multi-index $i \in \{0,1,2\}^d$, and $\lambda(i)$ is a scaling factor such that
\[ \lambda(i) = \begin{cases} 1 & \text{if $i = (0,\ldots,0)$,} \\ \bigl(\sum_{\ell = 1}^d i_\ell \bigr)^{-\alpha} & \text{otherwise.} \end{cases} \]
Here, $\alpha > 0$ is a decay factor controlling the smoothness of the polynomial $\phi$ via the damping of higher-order modes. With this set up, we define a class of randomly generated polynomials by constructing polynomials of the form \eqref{eq:phipoly}, with coefficients $c_i$ randomly and independently drawn from the uniform distribution on $[-1,1]$. In particular, we choose $\alpha = 2$, empirically chosen to represent a reasonable spread between simple and smooth interfacial geometry, to complex geometry involving high degrees of curvature.\footnote{In particular, $\alpha \geq 3$ tends to bias towards interfacial geometry with mild curvature, $\alpha \leq 1$ tends to solely create complex interfacial geometry, while $\alpha = 2$ ``evenly'' spreads between these characteristics.} We further characterise the class as follows:
\begin{itemize}
\item Let ${\mathcal A}$ denote the class of polynomials defined by \eqref{eq:phipoly}, with $\alpha = 2$ and $c_i$ uniformly randomly chosen from the interval $[-1,1]$, such that the following two conditions hold: (i) the volume fraction of $\{\phi < 0\}$ in $U = (-1,1)^d$ is between 10\% and 90\%; and (ii) using masking operations, it can be determined there exists a direction $e_k$ such that there are no branching points in the interface, i.e., $\phi$ and $\partial_k \phi$ do not share a real root in $U$.
\item Let ${\mathcal B}$ be defined in the same way as ${\mathcal A}$, except that condition (ii) cannot be verified, i.e., masking operations determine that each direction $e_k$ (likely) has a branching point in its corresponding implicitly defined height function, for every $1 \leq k \leq d$.
\end{itemize}
Thus, ${\mathcal A}$ and $\mathcal B$ represent two classes of randomly generated polynomials, differing according to the automated strategy of choosing a height direction discussed in \cref{sec:height}---for polynomials in class $\mathcal A$, there is an elimination axis for which it can be provably determined that the associated implicitly defined (multi-valued) height function is devoid of branching points, whereas polynomials in class $\mathcal B$ are likely to have vertical tangents for every possible height direction. As such, the two classes roughly separate the ``severity'' of the implicitly defined geometry in its representation as an implicitly defined height function. Meanwhile, the volume fraction condition is simply used to eliminate cases in which $U \cap \Gamma$ is empty or nearly empty.

\begin{figure}%
\centering\includegraphics{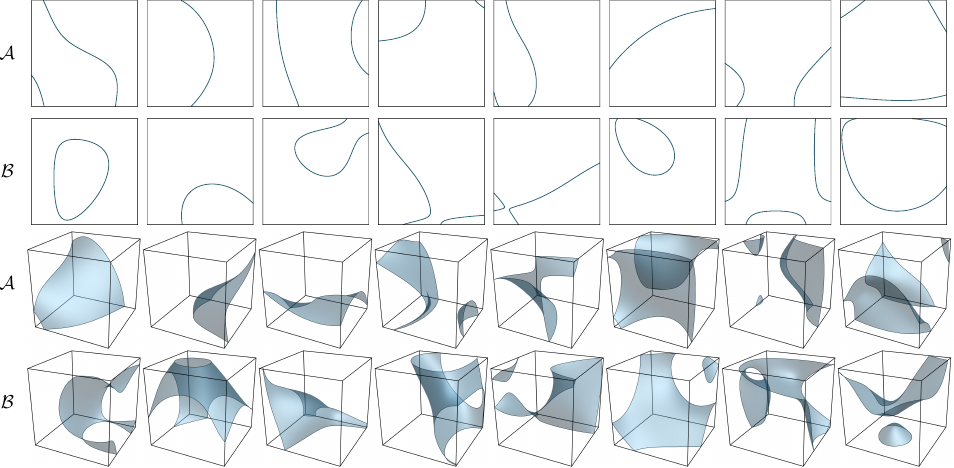}%
\caption{Examples of randomly generated interfacial geometry, corresponding to ``typical'' instances of the randomly generated polynomials of classes $\mathcal A$ and $\mathcal B$ defined in \cref{sec:random}. The top and bottom two rows correspond to $d = 2$ and $d=3$, respectively.}%
\label{fig:randomexamples}%
\end{figure}

Examples of the randomly generated polynomials are illustrated in \cref{fig:randomexamples}. One may observe that for each instance of class $\mathcal A$, there is a direction in which the corresponding height function has no vertical tangents, even if the height function is multi-valued. The examples of class $\mathcal B$ are generally more complex and contain vertical tangents for every possible elimination axis.

Next we define an error metric by which to assess the accuracy of the computed quadrature schemes on the above randomly generated geometries. Here we define a metric having close connection to cut-cell/embedded boundary finite element methods. Specifically, the quadrature algorithm is used to compute the mass matrix of a rectangular/prismatic element which has been cut by the implicitly defined geometry. Define the matrix-valued quantity $I^\pm$, whose $(i,j)$th entry is given by
\begin{equation} \label{eq:Iij} I_{ij}^{\pm} = \int_{U \cap \{\pm \phi < 0 \}} f_i f_j, \end{equation}
where the functions $f_i$ enumerate the polynomials of the Gauss-Lobatto nodal basis of biquadratic ($d = 2$) or triquadratic ($d = 3$) polynomials on $U = (-1,1)^d$; for details, see, e.g., \cite{ImplicitMeshPartOne}. Therefore, $I^+$ computes the mass matrix of the cut element on the positive side of the interface $\Gamma = \{\phi = 0\}$, and $I^-$ computes the mass matrix on the negative side. Let $\smash{I_q^\pm}$ denote the mass matrix computed by the quadrature scheme using a quadrature order $q$, and let $\smash{I_\text{ref}^\pm}$ denote a reference solution.\footnote{In the following results, the reference solution is computed using a large number of quadrature points, specifically $q = 100$ per sub-interval in \cref{algo:integrand}.} A relative error $E_q$ is defined such that
\begin{equation} \label{eq:Eq} E_q = \frac{ \| I_q^+ - I_\text{ref}^+\|_F + \| I_q^- - I_\text{ref}^- \|_F} { \|I_\text{ref}^+\|_F + \|I_\text{ref}^-\|_F}, \end{equation}
where $\| \cdot \|_F$ is the Frobenius norm. As such, $E_q$ measures the relative error in the computed mass matrices of both phases, $U \cap \{\phi < 0\}$ and $U \cap \{\phi > 0\}$.

The mass matrix error metric is used to examine the effectiveness of the quadrature algorithm on randomly generated geometry, as follows:
\begin{itemize}
\item One thousand instances\footnote{Experiments indicated that only a few hundred instances are needed to stabilise the computed statistics.}  of the polynomials in classes $\mathcal A$ and $\mathcal B$ are randomly generated (2000 total in 2D and 2000 in 3D). In every case, the automated strategy discussed in \cref{sec:height} is used to determine the height direction, at every level of the dimension reduction algorithm.
\item The kind of quadrature scheme to use in \cref{algo:integrand}, whether it is Gauss-Legendre or tanh-sinh, is defined \textit{a priori}, and shall depend only on the level of the dimension reduction process. Specifically, in the 2D case, the notation $(\mathsf{b},\mathsf{a})$ means that quadrature of type $\mathsf{b}$ is used on the outer integral, and $\mathsf a$ on the inner integral; in the 3D case, the notation $(\mathsf c, \mathsf b, \mathsf a)$ means quadrature of type $\mathsf c$ is used on the outer-most integral, $\mathsf a$ on the inner-most integral of the first eliminated axis, and $\mathsf b$ on the second eliminated axis. In particular, the objective here is to compute the entries of a mass matrix for which the integrand of the inner-most integral is the product of two polynomials, see \eqref{eq:Iij}; consequently, a natural choice is to apply Gauss-Legendre on the inner-most integral,\footnote{When $q \geq 3$, Gauss-Legendre quadrature will exactly integrate the inner-most vertical line integrals, the latter involving the product of two quadratic polynomials. Tanh-sinh, although rapidly convergent, will not be exact for such integrands.} and so $\mathsf a = \mathsf{GL}$ in the following tests.

\item For each polynomial and each combination of quadrature scheme type, the relative error $E_q$ is computed using quadrature orders ranging from $q = 1$ up to $q = 80$. In general, these quantities will decrease exponentially, up to the limits of machine precision, with a rate depending on the specific polynomial. To examine the corresponding distribution of convergence rates, we examine, as a function of $q$, the first, second, and third quartiles\footnote{Quartile metrics are better suited in this analysis than, say, the mean and variance; for example, the exponential slope of the third quartile is unaffected when the first quartile errors begin to plateau around machine precision.} of $E_q$, aggregated over all randomly generated polynomials in each particular class. %
\end{itemize}

\subsubsection{Two-dimensional results} 

\begin{figure}%
\centering\includegraphics{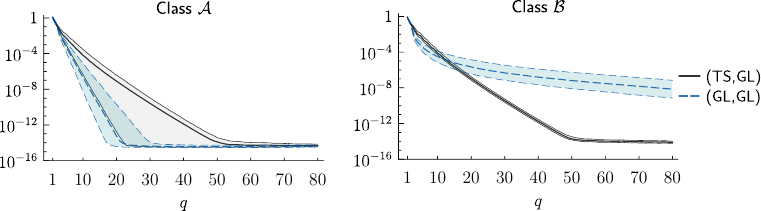}%
\caption{Quadrature accuracy as a function of $q$ for the two-dimensional randomly generated interfaces considered in \cref{sec:random}; some examples are shown in \cref{fig:randomexamples}. The quadrature combination $(\mathsf{GL},\mathsf{GL})$ is indicated in blue with dashed lines, whereas the combination $(\mathsf{TS},\mathsf{GL})$ is indicated in grey with solid lines. In each case, the shaded region demarcates the boundaries of the first and third quartiles of $E_q$, defined by \eqref{eq:Eq}, aggregated over all 1000 polynomial instances, while the interior line indicates the median.}
\label{fig:random2d}
\end{figure}

In 2D, we consider two combinations of quadrature schemes, $(\mathsf{GL},\mathsf{GL})$ and $(\mathsf{TS},\mathsf{GL})$. %
\cref{fig:random2d} illustrates the results of the study, divided according to the classes $\mathcal A$ and $\mathcal B$. Several observations can be made:
\begin{itemize}
\item For the randomly generated geometry of class $\mathcal A$, the quadrature accuracy using $(\mathsf{GL},\mathsf{GL})$ is generally far superior to $(\mathsf{TS},\mathsf{GL})$. Indeed, the corresponding interfacial geometry is the graph of a smooth height function devoid of branches, and as such we intuitively expect to benefit from the rapid convergence of high-order Gaussian quadrature. Moreover, when Gaussian quadrature functions effectively, it can be expected to yield more rapid convergence than tanh-sinh quadrature.
\item For polynomials in class $\mathcal B$, the quadrature choice $(\mathsf{GL},\mathsf{GL})$ performs very poorly, and does not reach machine precision. This should also be intuitively expected, as Gaussian quadrature performs poorly on integrands with endpoint singularities; in this case, the endpoint singularities occur in the derivative of the integrands, corresponding to vertical tangents of the height function. On the other hand, $(\mathsf{TS},\mathsf{GL})$ is effective at treating these kinds of singularities. However, note that for small $q$, around ten or less, $(\mathsf{GL},\mathsf{GL})$ can be substantially more accurate than $(\mathsf{TS},\mathsf{GL})$; in a sense, tanh-sinh quadrature can require a moderate number of quadrature points before its convergence accuracy is ultimately superior to Gaussian quadrature.
\item One can also observe that the spread in errors for $(\mathsf{TS},\mathsf{GL})$ on class $\mathcal B$ is quite small;  evidently, this combination of quadrature scheme behaves consistently across the wide range of possible interfacial geometry in class $\mathcal B$. 
\end{itemize}

\subsubsection{Three-dimensional results}

\begin{figure}%
\centering\includegraphics{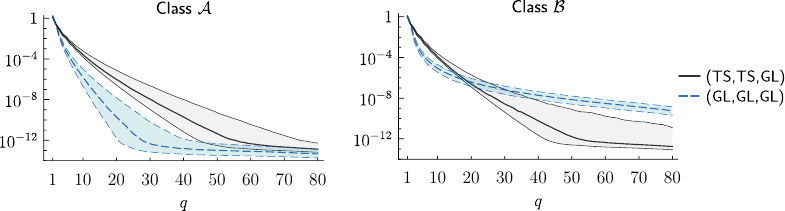}%
\caption{Quadrature accuracy as a function of $q$ for the three-dimensional randomly generated interfaces considered in \cref{sec:random}; some examples are shown in \cref{fig:randomexamples}. The quadrature combination $(\mathsf{GL},\mathsf{GL},\mathsf{GL})$ is indicated in blue with dashed lines, whereas the combination $(\mathsf{TS},\mathsf{TS},\mathsf{GL})$ is indicated in grey with solid lines. In each case, the shaded region demarcates the boundaries of the first and third quartiles of $E_q$, defined by \eqref{eq:Eq}, aggregated over all 1000 polynomial instances, while the interior line indicates the median.}
\label{fig:random3d}
\end{figure}

In 3D, there are four main combinations of quadrature schemes,  $(\mathsf{GL},\mathsf{GL},\mathsf{GL})$, $(\mathsf{TS},\mathsf{GL},\mathsf{GL})$, $(\mathsf{GL},\mathsf{TS},\mathsf{GL})$, and $(\mathsf{TS},\mathsf{TS},\mathsf{GL})$. Experiments indicate that the first two combinations yield similar results as each other, and likewise for the last two. As such, for brevity, the presented results are of $(\mathsf{GL},\mathsf{GL},\mathsf{GL})$ and $(\mathsf{TS},\mathsf{TS},\mathsf{GL})$, and are shown in \cref{fig:random3d}. Similar observations as the two-dimensional case can be made:
\begin{itemize}
\item For the randomly generated geometry of class $\mathcal A$, the quadrature accuracy using $(\mathsf{GL},\mathsf{GL},\mathsf{GL})$ is generally superior to $(\mathsf{TS},\mathsf{TS},\mathsf{GL})$. From the perspective of the first eliminated coordinate axis, the interfacial geometry of polynomials in $\mathcal A$ can be recast as the graph of a smooth height function devoid of branches, and as such we intuitively expect Gaussian quadrature to be effective. %
\item Similar to the 2D case, for polynomials in class $\mathcal B$, $(\mathsf{GL},\mathsf{GL},\mathsf{GL})$ performs very poorly, which is to be expected given that every coordinate direction yields a height function with (near) vertical tangents, i.e., endpoint singularities. In this setting, tanh-sinh quadrature is better equipped at handling arbitrary endpoint singularities, and can be seen to ultimately attain higher convergence rates. Though, similar to the 2D case, Gaussian quadrature often yields better accuracy when the number of quadrature points is small; in this case, the transition point is around $q \approx 20$. 
\item One may also observe the wider spread in convergence rates, particularly of class $\mathcal B$ polynomials, compared to the 2D case. This is reflective of the increased interfacial complexity in 3D (as all examples in \cref{fig:randomexamples} show), as well as the more difficult geometries obtained in the dimension-reduced base integrals. We note here that if the randomly generated geometry is made smoother, by increasing the exponent $\alpha$ in the damping of higher-order Legendre modes (see \cref{eq:phipoly}), then the convergence rates in both classes increase and the interquartile spread decreases; indeed, for $\alpha$ sufficiently large the randomly generated polynomials are almost linear and the implicitly defined domains essentially become polyhedral.
\end{itemize}

\subsubsection{Summary}

The above results, making use of randomly generated interfacial geometry, suggest the following guidelines for automatically choosing the type of one-dimensional quadrature in \cref{algo:main,algo:integrand}:
\begin{enumerate}[label=(\roman*),leftmargin=3em]
\item If it can be proven (using masking operations) that the polynomial(s) involved are the graph of a (potentially multi-valued) height function devoid of branching points/vertical tangents, then use Gauss-Legendre quadrature in the base integral; 
\item Otherwise, and if high precision is required, consider applying tanh-sinh quadrature;
\item However, if only modest accuracy is required or a relatively small quadrature order $q$ is to be used, then apply Gauss-Legendre quadrature.  
\end{enumerate}
Situations (i) and (iii) are applicable to the examples in class $\mathcal A$ as well as, for example, the vast majority of scenarios arising in cut-cell finite element methods of fluid-interface dynamics, see, e.g., \cite{InterfacialGauge,ImplicitMeshPartTwo}.

The results also show that the quadrature algorithm is effective at handling complex interfacial topology, such as multi-component domains and varying degrees of curvature, without having to resort to heuristic subdivision methods. Through a suitable combination of one-dimensional quadrature types, approximately exponential convergence is obtained, i.e., doubling $q$ approximately doubles the number of correct digits. The rate of convergence is related to the degree of curvature, and in the majority of cases, reasonably fast slopes are attained. Slower rates of convergence can be obtained in 3D, e.g., the results of \cref{fig:random3d}(right) exhibit a wider spread of convergence rates above the median slope; these are associated with a high degree of localised interface curvature, as discussed further in later results. In these more rare scenarios, a small amount of constraint cell subdivision (i.e., localised $h$-refinement) can be effective in restoring rapid convergence; see, e.g., the supplementary material which contains a numerical study demonstrating this. Constraint cell refinement opens up the possibility of using different height directions on different parts of the interfacial geometry. However, we recall that an overarching design principle of these quadrature schemes was to avoid a reliance on such methods.

\subsection{Bilinear example}
\label{sec:bilinear}

\begin{figure}%
\centering\includegraphics{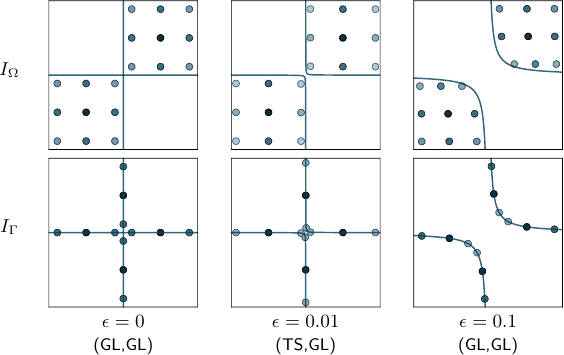}%
\caption{Illustrations of the $q = 3$ quadrature schemes for the bilinear test problems considered in \cref{sec:bilinear}; in each instance, the quadrature points are shaded according to their weight, such that low-valued weights are pale and higher-valued weights are darker.}
\label{fig:bilinearquads}
\end{figure}

\begin{figure}%
\centering\includegraphics{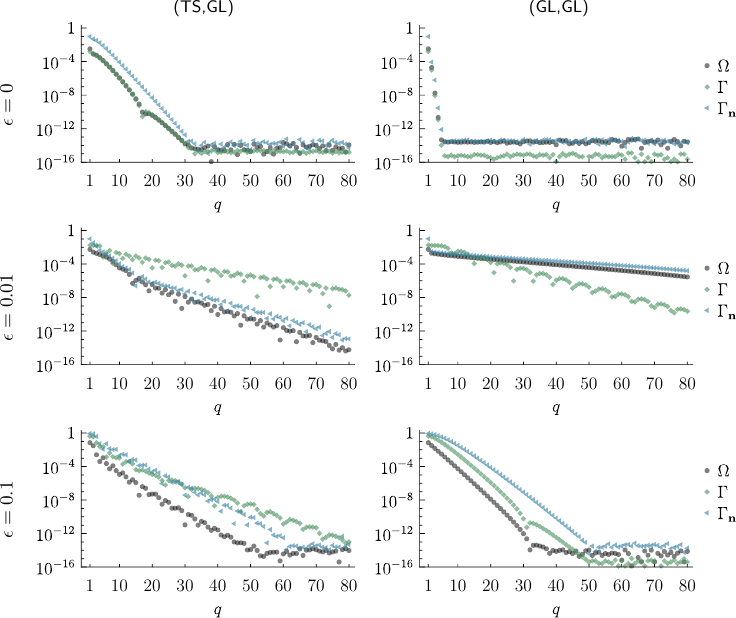}%
\caption{Quadrature accuracy as a function of $q$ for the bilinear test problems considered in \cref{sec:bilinear}. Here, $\bullet$ (black) denotes the relative error in the computed volumetric integral, $I_\Omega$; $\squaresymbol$ (green) denotes the relative error in the computed surface integral in non-flux form, $I_\Gamma$; and $\trianglesymbol$ (blue) denotes the relative error in the computed flux-form surface integral, $I_{\Gamma_{\textbf n}}$.}%
\label{fig:bilineartests}%
\end{figure}

In the next example, we consider a simple bilinear polynomial to study the effect of localised, high degree curvature. Consider the polynomial $\phi(x,y) = (x-\tfrac12)(y-\tfrac12) - \epsilon^2$, whose interface depends on $\epsilon$ as follows: for $\epsilon$ sufficiently large, $\Gamma$ exhibits a smooth rounded corner; as $\epsilon \to 0^+$, a sharp corner of curvature ${\mathcal O}(\epsilon^{-1})$  develops; and for $\epsilon = 0$, $\Gamma$ is a cross centred at $(0.5,0.5)$. Similar to the ellipsoidal test problem, we consider three kinds of integral quantities,
\[ I_{\Omega} := \int_{U \cap \{\phi > 0\}} f, \quad I_{\Gamma} := \int_{U \cap \{\phi = 0\}} f, \quad I_{\Gamma_{\textbf n}} := \int_{U \cap \{\phi = 0\}} f \mathbf n, \]
where $U := (0,1)^d$ is fixed, and the integrand $f : \R^2 \to \R$ is given by $f(x) = \cos (\tfrac14 \|x\|_2^2)$. We consider three values of $\epsilon$, specifically $0$, $0.01$ and $0.1$, representing the degenerate case, high localised curvature, and smooth geometry, respectively. \cref{fig:bilinearquads} illustrates the corresponding quadrature schemes for $q = 3$, while \cref{fig:bilineartests} contains the results of a $q$-refinement study for the three integral quantities, using two kinds of quadrature choices, $(\mathsf{TS},\mathsf{GL})$ and $(\mathsf{GL},\mathsf{GL})$ (here we use the notation as in \cref{sec:random}). Several observations can be made:
\begin{itemize}
\item When $\epsilon = 0$, the quadrature algorithm correctly handles the singular interface geometry and produces an intuitive quadrature scheme: for the volumetric integral, a tensor-product scheme is produced for the upper-right and lower-left quadrants, and for the surface integral, the four arms of $\Gamma$ are sampled with a simple one-dimensional scheme; see \cref{fig:bilinearquads}, left.
\item For the relatively smooth interface corresponding to $\epsilon = 0.1$, the fastest convergence accuracy is obtained when Gauss-Legendre quadrature is used in all directions, as intuitively expected.
\item In the high-curvature case of $\epsilon = 0.01$, exhibiting a rounded corner of radius about 1\% of the cell width, the convergence rates of both schemes slow. In essence, there is a thin transition layer of width $\mathcal O(\epsilon)$ in the behaviour of the associated height functions. In the case of Gaussian quadrature, i.e., $(\mathsf{GL},\mathsf{GL})$, this severely impacts its convergence rate. However, tanh-sinh quadrature is well-equipped at handling almost arbitrary endpoint singularities, helping to explain why the convergence rates of the $(\mathsf{TS},\mathsf{GL})$ scheme are not as severely impacted.
\item In the case of $\epsilon = 0.01$ and the $(\mathsf{TS},\mathsf{GL})$ scheme, note also the transition between a steep slope (mainly of the volumetric and surface-flux integrals) prior to $q \approx 15$ to a shallower slope after $q \approx 15$. This transition relates to the convergence of the interfacial geometry to the degenerate cross shape as $\epsilon \to 0$. For smaller and smaller $\epsilon$, more and more quadrature points are needed to detect the diminishingly-small rounded corners of the interface; on the other hand, the integral error incurred by ignoring this geometry also diminishes as $\epsilon \to 0$. As such, the transition point ($q \approx 15$ for $\epsilon = 0.01$) will move to larger values of $q$ as $\epsilon$ is made smaller, restoring the rapid convergence and effectively obtaining the results shown in \cref{fig:bilineartests}, top row, once $\epsilon$ is sufficiently small.
\end{itemize}
Combined, these observations indicate that an appropriate choice of tanh-sinh and Gauss-Legendre quadrature can robustly handle a wide degree interface curvature. The above analysis considered these effects in the case the quadrature algorithm is ``blind'' to the interfacial geometry at hand. Naturally, if the geometry is known ahead of time, then one can subdivide appropriately so as to apply adaptive quadrature methods to achieve higher rates of convergence. Further comments about potential automated strategies for doing so are given in \cref{sec:conclusion}.

\subsection{Trilinear tunnel}
\label{sec:trilinear}

\begin{figure}%
\centering\includegraphics{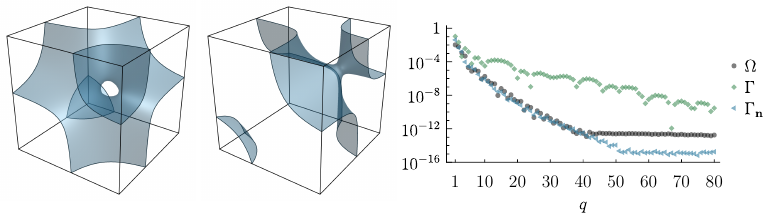}%
\caption{\textit{(left pair)} Two different views of the interfacial geometry defined by the trilinear polynomial of \cref{sec:trilinear}; note the tunnel connecting one vertex of the unit cube $(0,1)^3$ to the interior. \textit{(right)} Quadrature accuracy as a function of $q$ for the trilinear tunnel problem considered in \cref{sec:trilinear}. Here, $\bullet$ (black) denotes the relative error in the computed volumetric integral; $\squaresymbol$ (green) denotes the relative error in the computed surface integral in non-flux form; and $\trianglesymbol$ (blue) denotes the relative error in the computed flux-form surface integral.}
\label{fig:trilinear}
\end{figure}

The next example demonstrates that even a low degree polynomial can yield non-trivial interfacial topology. Specifically, consider the trilinear polynomial $\phi: \R^3 \to \R$ given by 
\[ \phi(x,y,z) = 0.5 - 1.4 z + 2.9 x y - 6.5 x y z + 3.2 x z - 1.2 x + 3.3 y z - 1.3 y. \]
The interface of $\phi$ in the unit cube $(0,1)^3$ has two components, one of which exhibits a tunnel, see \cref{fig:trilinear}. Applied to the integrand $f(x) =  \cos (\tfrac14 \|x\|_2^2)$, \cref{fig:trilinear} also illustrates the results of a $q$-refinement study similar to those in prior studies, i.e., of the volumetric integral, surface integral in non-flux form, and flux-form surface integral of $f$. In particular, following the guidelines developed in \cref{sec:random}, the computed results apply $(\mathsf{TS},\mathsf{TS},\mathsf{GL})$ quadrature, i.e., Gauss-Legendre on the first eliminated axis, and tanh-sinh quadrature on every base integral. The results show similar convergence rates as in prior examples---once again, we observe that surface integrals computed in flux-form can be far more accurate than surface integrals in non-flux form.

\subsection{Surfaces with singularities}
\label{sec:singular}

\begin{figure}%
\centering\includegraphics{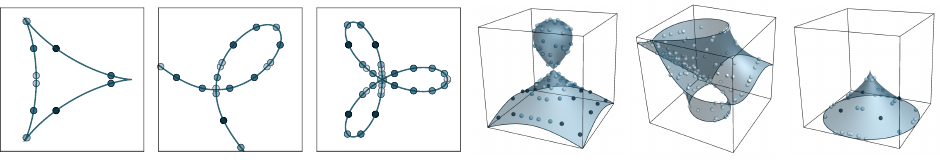}%
\caption{Examples of the singularity test problems included in the supplementary material.}%
\label{fig:smsingular}%
\end{figure}

The supplementary material accompanying this paper provides additional examples involving singularities, such as self-intersections and cusps; some are shown in \cref{fig:smsingular}. These examples are more contrived in the sense that they typically would not occur in computational physics problems, say; nonetheless, they serve to stress-test the quadrature algorithm in situations involving roots of high-order multiplicity.

\subsection{Application to simplices and intersecting geometry}

By design, the quadrature algorithm can take in as input multiple multivariate polynomials. In this scenario, the implicitly defined interface is the union of each polynomial's zero isosurface, and the output of the algorithm is a quadrature scheme that can be used to compute integrals on individual connected components of the corresponding implicitly defined domain. We demonstrate this utility here, through some examples involving simplices and multi-component domains of intersecting geometry.

\subsubsection{Extension to simplices}
\label{sec:simplices}

\begin{figure}%
\centering\includegraphics{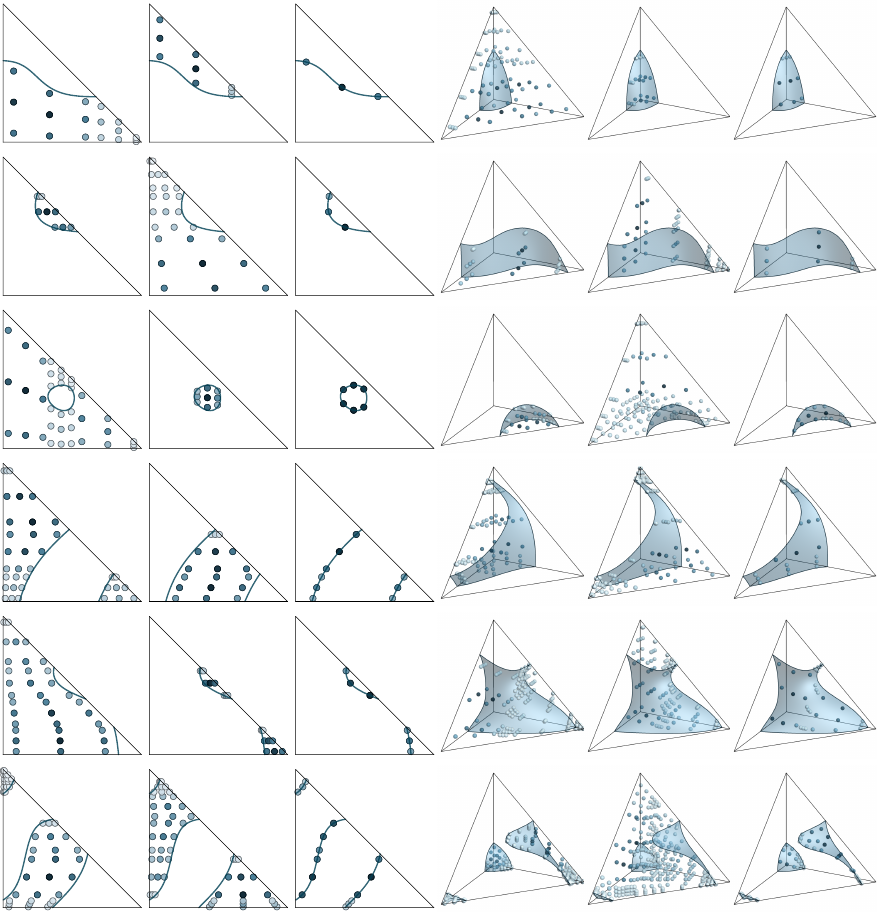}%
\caption{Examples illustrating the adaptation of the quadrature algorithm to simplex constraint cells, following the mechanism discussed in \cref{sec:simplices}. In each group of three, the left and middle figure corresponds to volumetric integrals over the implicitly defined regions $V \cap \{\phi < 0\}$ and $V \cap \{\phi > 0\}$, whereas the right figure corresponds to the surface integral over $V \cap \{\phi = 0\}$, where $V$ is the standard triangle or tetrahedron simplex; in each instance, the quadrature points are shaded according to their weight, such that low-valued weights are pale and higher-valued weights are darker.}%
\label{fig:simplices}%
\end{figure}

The dimension reduction process operates most naturally when the constraint cell $U \subset \R^d$ is a hyperrectangle. To handle the case of a simplex constraint cell, a simple approach is to extend the simplex to its enclosing hyperrectangle and then incorporate an additional polynomial $\psi : \R^d \to \R$ whose purpose is to cut the implicitly defined geometry along the sloping hyperplane of the target simplex. For example, suppose we wish to compute a quadrature scheme for
\[ \int_{V \setminus \{\phi = 0\}} f, \]
where $V := \{ x \in \R^d : x_i \geq 0 \text{ and } \sum_i x_i \leq 1 \}$ is the standard unit simplex and $\phi$ is a polynomial. Define the linear polynomial $\psi$ by $\psi(x) = -1 + \sum_i x_i$, and execute\footnote{It is also worthwhile to modify the masks of all input polynomials such that the masks are disabled, i.e., value $0$, in all subcells outside $V$; doing so helps to eliminate unneeded zeros of $\mathcal P$, $\mathcal Q$, or $\mathcal R$, i.e., helps to ignore the irrelevant branching or intersection points of the multi-valued height function which are outside the simplex.} \cref{algo:main} to compute a volumetric quadrature scheme for the domain $U \setminus (\{\phi = 0 \} \cup \{\psi = 0\})$, where $U = (0,1)^d$. This yields a volumetric quadrature scheme filling all of $U$; one may then simply discard any quadrature nodes outside the simplex and leave the quadrature weights of all remaining nodes unmodified. This construction yields a high-order quadrature scheme appropriate for computing integrals over the connected components of $V \cap \{\phi \neq 0\}$. Examples are shown in \cref{fig:simplices}. We note that \cref{algo:main,algo:integrand} could easily be modified to directly incorporate $\psi$ for increased efficiency, e.g., so that the quadrature points in $U \setminus V$ are not needlessly computed.

\subsubsection{Intersecting geometry}
\label{sec:intersections}

\begin{figure}%
\centering\includegraphics{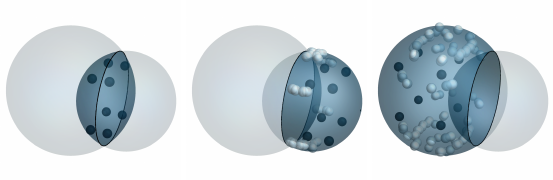}%
\caption{An example showing how the quadrature algorithm can be used to build quadrature schemes for boolean operations on implicitly defined geometry, in this case of two overlapping spheres. In each figure, a volumetric quadrature is shown (using $q = 2$), for the lens-shaped intersection of both balls, and the two crescent-shaped regions (where one ball is subtracted from the other); the quadrature points are shaded according to their weight, such that low-valued weights are pale and higher-valued weights are darker.}
\label{fig:circles}
\end{figure}

\begin{figure}%
\centering\includegraphics{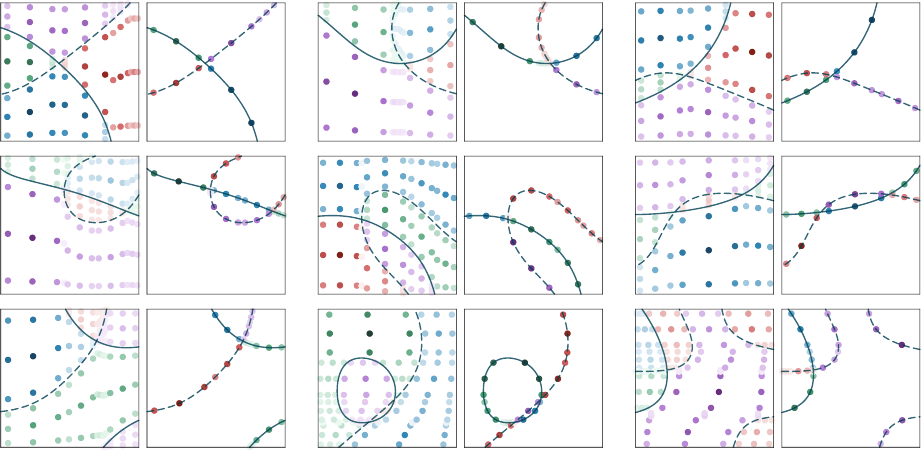}%
\caption{Examples illustrating how the quadrature algorithm can be used to build quadrature schemes of intersecting geometry. In each pair, the left figure corresponds to volumetric integrals while the right figure corresponds to surface integrals. In each instance, the quadrature points are coloured according to the signs of the input polynomials; each cluster of identical-coloured quadrature points yields a high-order quadrature scheme of the corresponding region.}
\label{fig:intersection}
\end{figure}

Finally, we show how the quadrature algorithm can be used to compute quadrature schemes for intersecting or overlapping domains. A general purpose mechanism is to: (a) invoke \cref{algo:main} on the set of polynomials implicitly defining each domain, and then: (b) cluster the resulting quadrature points according to the signs of the input polynomials. This mechanism can be used to compute quadrature schemes for the various boolean operations, \textsc{and}, \textsc{or}, \textsc{xor}, etc. By construction, each cluster yields a high-order accurate quadrature scheme when applied to sufficiently smooth integrands, yielding similar kinds of $h$-refinement and $q$-refinement convergence as shown in the other numerical experiments of this work. To illustrate the kind of schemes produces, \cref{fig:circles} demonstrates an example involving two overlapping spheres and shows the volumetric quadrature schemes for the crescent- and lens-shaped regions; a similar clustering approach can be used to compute surface quadrature rules for these regions. Additional examples are shown in \cref{fig:intersection} of intersecting geometry implicitly defined by randomly generated 2D polynomials.\footnote{The same experiment can be done in 3D but the resulting geometry can be highly intricate, the visualisation of which is better suited for interactive plots rather than what can be shown on paper.}

\section{Concluding Remarks}
\label{sec:conclusion}

The quadrature algorithms developed in this work reinterpret the implicitly defined geometry as the graph of an implicitly defined, multi-valued height function. Using the corresponding splitting of coordinates, the target integral is reformulated as a base integral, in one fewer dimensions, whose associated integrand $\mathcal F$ is either a vertical line integral in the height direction (in the case of volumetric integrals), or a point-wise interfacial evaluation of the user-supplied integrand, taking into account curved surface Jacobian factors (in the case of surface integrals). Continuing recursively, the dimension reduction process eliminates one coordinate axis at a time, down to a one-dimensional base case; the base case then performs a one-dimensional quadrature, which in turn initiates a recursion back up the tree, with each new integrand evaluation performing a one-dimensional quadrature. To achieve a high order of accuracy, it is important to appropriately identify and treat the piecewise smooth behaviour of these one-dimensional integrals. This can be achieved through two mechanisms: (i) by dividing the domain of the base integral into one or more pieces, each one corresponding to the different characteristics of the multi-valued height function (e.g., crossing and branching points); and (ii) by appropriately choosing a quadrature scheme for the one-dimensional integrals (e.g., tanh-sinh quadrature to effectively handle endpoint singularities associated with branching points). In essence, the output of the algorithm is a quadrature scheme consisting of one or more patches, with each patch containing a tensor-product quadrature scheme warping to the implicitly defined geometry.\footnote{We note here that the dimension reduction process is in some ways similar to the algorithm pioneered by Collins \cite{10.1007/3-540-07407-4_17} for computing the \textit{cylindrical algebraic decomposition} (CAD), a well-known tool in the field of real algebraic geometry \cite{caviness2012quantifier,BasuPollackRoy}. Roughly speaking, a CAD is a decomposition of $\R^d$ into cells, such that a given set of polynomials have constant sign throughout each cell. A typical algorithm for computing a CAD operates through quantifier elimination and a projection phase (similar to finding a height direction and computing resultants in the description of the base integral) down to a one-dimensional base case, combined with a lifting phase which recursively adds dimensions to finish with a decomposition of the full space $\R^d$. Moreover, under certain conditions, each cell of a CAD is diffeomorphic to an open hypercube $(0,1)^\ell$, for some $\ell \leq d$. The patches mentioned above may be viewed as the image, under the action of the diffeomorphism, of a canonical tensor-product scheme in $(0,1)^\ell$, suitably modified to account for the (implicitly defined) Jacobian determinant of this mapping. We also note that, in the present setting, we are considering implicitly defined geometry in a user-defined constraint cell $U$; a CAD could do the same if its input polynomials included a description of the boundary of $U$. For example, in 2D, if $U = (x_1,x_2) \times (y_1,y_2)$, then append the polynomials $x - x_1$, $x - x_2$, $y - y_1$, and $y - y_2$. In a sense, the quadrature algorithm directly incorporates this geometry, rather than have it be specified through additional polynomials. Finally, we note that the quadrature algorithm restricts attention to $U$, thus removing the computational expense of computing a full CAD, while the use of masking operations helps to eliminate additional CAD cells outside $U$.} 

Key findings include the following:
\begin{itemize}
\item Under the action of $h$-refinement and with $q$ fixed, for sufficiently smooth geometry and integrands (e.g., away from singularities such as cusps), the order of accuracy of the quadrature scheme is approximately $2q$.
\item Under the action of $q$-refinement, i.e., keeping both the constraint cell $U$ and input polynomials $\{\phi_i\}$ fixed, the convergence is (asymptotically) approximately exponential in $q$, i.e., doubling $q$ approximately doubles the number of accurate digits. In particular, the convergence behaviour and rate depends on the degree of curvature in the interface: smooth interfacial geometry yields rapid convergence (similar to Gaussian quadrature), whereas localised regions of very high curvature (e.g., nearly-sharp corners) can slow convergence; see, e.g., \cref{sec:bilinear}.
\item Surface quadrature schemes in flux form, i.e., to compute $\int_\Gamma f \mathbf n$, are generally more accurate than those in non-flux form, i.e., $\int_\Gamma f$, relating to the property that the additional factor of $\mathbf n$ can remove (or cancel entirely) the behaviour of curved surface Jacobians. Fortunately, in many applications of surface quadrature, e.g., discontinuous Galerkin methods, it is more important to achieve high levels of accuracy in flux form than it is in non-flux form.
\item In more contrived examples involving interfacial singularities, especially cusps, the numerical conditioning of the quadrature problem worsens, and this can impact the maximal possible accuracy in finite precision arithmetic. For example, the results of \cref{sec:singular} show that some integrals can be limited to a maximum accuracy of ${\mathcal O}(\epsilon^\alpha)$, where $\epsilon$ is machine epsilon and $\alpha < 1$. This limited precision is related to the conditioning of computing roots near the singularities; at the singularities themselves, the roots are typically of high multiplicity. Here, two aspects are worth recalling: from a well-posedness standpoint, such problems are highly sensitive to perturbations in the input polynomials, e.g., small perturbations in a cusp can result in large changes in its surface area. On the other hand, these aspects of ill-conditioning do not always impact integral accuracy, e.g., full machine precision can be obtained for problems involving corners and junctions arising from crossing interfaces.
\end{itemize}

In this work we have mainly explored the use of Gauss-Legendre and tanh-sinh quadrature. In general, Gauss-Legendre performs well when the associated height functions do not have a vertical tangent, and, even if they do, it often yields a comparatively high accuracy when $q$ is small. On the other hand, tanh-sinh is better equipped at handling end-point singularities, and is therefore more effective for height functions exhibiting vertical tangents, such as at branching points. Numerical experiments using randomly generated interfacial geometry, see \cref{sec:random}, support this assessment and suggest the following strategy for choosing between Gauss-Legendre and tanh-sinh: if masking operations indicate there are branching points (i.e., the set $\mathcal P$ discussed in \cref{sec:driver} is likely non-empty), and high precision is sought, then apply tanh-sinh quadrature in the base integral; otherwise, or if the target application is limited to a small $q$, use Gauss-Legendre quadrature.

From the outset, the quadrature algorithms were designed so as to not rely on ad hoc subdivision methods. Nevertheless, in some applications it could be highly beneficial to subdivide the constraint cell $U$ into smaller subcells (of the same shape). Doing so effectively allows the choice of height direction/elimination axis to become more locally adaptive to the interfacial geometry at hand, such that different portions of the geometry use different height directions. For example, in the vast majority of cases of the randomly generated geometry in \cref{sec:random}, just one or two levels of subdividing $U$ into equal sized quadrants/octants, say, would yield the kinds of geometry for which Gauss-Legendre yields rapid convergence, faster than what is possible without any subdivision; see, e.g., the numerical results provided in the supplementary material. Of course, subdivision may incur a greater number of quadrature points. However, it should be noted that one cannot always rely on this kind of subdivision in order to resolve complex geometry of $\Gamma$; for example, singularities require a direct treatment and cannot be efficiently handled solely through constraint cell subdivision.

\begin{table}%
\centering\footnotesize%
\begin{tabular}{lcccccc} \hline
\multirow{2}{*}{Test problem} & \multirow{2}{*}{Build} &  \multicolumn{5}{c}{Order $q$} \\
 &  & 1 & 2 & 4 & 8 & 16 \\ \hline
2D ellipse, $\mu\text{s}$ per grid cell with $U \cap \Gamma \neq \emptyset$ & 0.2 & 0.1 & 0.2 & 0.5 & 1.0 & 2.5 \\
3D ellipsoid, $\mu\text{s}$ per grid cell with $U \cap \Gamma \neq \emptyset$ & 1.5 & 0.2 & 0.8 & 2.2 & 8.8 & 50 \\
2D randomly generated, average $\mu\text{s}$ per instance  & 0.4 & 0.2 & 0.3 & 0.6 & 1.2 & 2.9 \\
3D randomly generated, average $\mu\text{s}$ per instance & 6 & 0.9 & 2.2 & 8.0 & 25 & 120 \\ \hline
\end{tabular}
\caption{Timing measurements, in microseconds, of the build phase (i.e., execution of \cref{algo:main}), and construction phase (i.e., execution of \cref{algo:integrand}) for the indicated studies; see the main text for detail, \cref{sec:conclusion}. Experimental results obtained on an Intel Xeon E3-1535m v6 laptop, single core, operating at approximately 4Ghz.}%
\label{tab:timing}\vspace{-2em}%
\end{table}

We have yet to discuss the computational cost of the quadrature algorithm; some comments and timing experiments are provided here. In short, the quadrature algorithm is fast for input polynomials of modest degree and can significantly outperform the algorithm developed in the author's prior work \cite{HighorderImplicitQuad,algoim}. To provide an indication of the computational cost, \cref{tab:timing} measures the time, in microseconds, needed to construct and evaluate the quadrature scheme for different values of $q$. In particular, the ``Build'' column measures the time needed to execute \cref{algo:main} (computing masks, face restrictions, resultants, pseudo-discriminants, etc.) whereas the other columns indicate the time required to evaluate \cref{algo:integrand} (evaluating vertical line restrictions of multivariate Bernstein polynomials, computing roots, mask filtering, sorting into subintervals, etc.), assuming that the cost of evaluating the inner-most user-supplied integrand is negligible. The first two rows correspond to the $h$-refinement study of the ellipsoid in \cref{sec:ellipse}, and measures the time per grid cell, in microseconds, needed to compute volumetric quadrature schemes, averaged over all grid cells containing the interface; meanwhile, the last two rows correspond to the randomly generated examples in \cref{sec:random} and measures the time per instance, in microseconds, needed to compute volumetric quadrature schemes, averaged over all randomly generated instances. Although the number of quadrature points scales as ${\mathcal O}(q^d)$, the timings in \cref{tab:timing} do not directly follow this trend. This is because the quadrature formation is spread across a variety of components, such as calculating polynomial roots of mixed degree polynomials (sometimes using a fast hybrid of Bernstein subdivision plus Newton's method, sometimes using more complex eigenvalue methods), mask filtering operations, root sorting, and formation of the Gauss-Legendre/tanh-sinh quadrature points, all of which have different asymptotic complexities. Generally speaking, profiling experiments indicate that (at least for the test problems under consideration) polynomial root finding is not a bottleneck; instead, most of the time is spread between multivariate Bernstein polynomial manipulation, masking operations, and formation of the tensor-product scheme. We also note that our C++ implementation is not particularly optimised; nevertheless, the timing measurements reported in \cref{tab:timing} are at least as fast as the identical ellipsoidal test carried out in \cite{HighorderImplicitQuad}, and in some cases more than three times faster.

We have not carefully considered here the issue of exploding polynomial degrees arising from the use of resultant and discriminant methods. In particular, the resultant of two tensor-product polynomials of degree $(p,\ldots,p) \in \mathbb N^d$ yields a tensor-product polynomial of maximal degree $(2p^2, \ldots, 2p^2) \in \mathbb N^{d-1}$ (see \cref{subsec:resultants}). Therefore, in a worst case scenario, a 2D quadrature problem, with input two polynomials of degree $(p,p)$, may result in a degree $2p^2$ polynomial implicitly defining the domain of its base integral. In 3D, this effect compounds: a quadrature problem with input two polynomials of degree $(p,p,p)$ may result in a 1D base integral (after elimination of two coordinate axes) involving a polynomial of degree $2 (2 p^2)^2 = 8 p^4$. Meanwhile, the computational cost of building these polynomials has asymptotic complexity at least the degree squared, and so the worst-case complexity of the quadrature algorithm is at least $\mathcal O (p^4)$ in 2D and $\mathcal O(p^8)$ in 3D. (The preceding analysis ignores some aspects, such as perhaps needing eigenvalue methods to compute the roots of these polynomials, potentially with higher complexity.) Evidently, this cost, as well as associated aspects of numerical conditioning, may limit the practicality of the dimension reduction algorithm when the input polynomials describing the implicitly defined geometry are of high degree. Consequently, the quadrature algorithm developed in this work is perhaps best suited to applications involving low-to-medium degree multivariate polynomials. In any case, numerical experiments indicate that the algorithm is quite robust with respect to manipulating high degree polynomials, owing in part to the use of the Bernstein basis; see, e.g., some of the test problems included in the supplementary material.

Finally, we mention here some possibilities of further research that could lead to more advanced implementations of these quadrature algorithms:
\begin{itemize}
\item Several promising approaches could be implemented to increase quadrature adaptivity. In this work, and mainly for simplicity, we used a fixed number of quadrature points $q$ per one-dimensional subinterval (see \cref{sec:1dquad}). Clearly, this could lead to an unnecessary number of quadrature points in the case of small subintervals. Thus, one could develop more adaptive methods that expend fewer points in small subintervals, leading to a greater control in bounding the total number of output quadrature points. Extending this idea further, a more advanced implementation could begin with a low-order scheme, and adaptively add quadrature points where needed based on an error indicator, in the spirit of methods like Gauss-Kronrod, or through the re-use of existing points as afforded by, e.g., the nesting already present in tanh-sinh schemes. Guiding these efforts, one could use estimates of radius of convergence of the analyticity of the implicitly defined height functions, as quantified by the implicit function theorem. All of these latter ideas can exploit the fact the setup phase of the algorithm, i.e., \cref{algo:main}, is independent of the scheme evaluation phase, i.e., \cref{algo:integrand}.
\item We explored here the use of tanh-sinh quadrature, owing to its effectiveness in treating almost arbitrary endpoint singularities. However, the implicitly defined height functions, here associated with domains implicitly defined by multivariate polynomials, have additional structure that could be exploited. For example, the endpoint singularities are typically associated with integrands which algebraically decay to zero. As such, alternatives to tanh-sinh quadrature which specifically target the behaviour of multi-valued height functions may yield higher rates of convergence.\footnote{Gauss-Jacobi is one possibility, but this approach requires precise determination of the algebraic power of the end-point singularity, which could vary from segment to segment of the height function, thus making the approach intricate to implement.}
\item The advances made with the above ideas could also open up other kinds of strategies for choosing the height direction/elimination axis, beyond the simple (though often effective) methods outlined in \cref{sec:height}.
\item A relatively straightforward method has been used to compute resultant and pseudo-discriminant polynomials (see \cref{sec:resultants}), whose roots are then computed for use in the base integral. It may be possible to avoid the explicit construction of these polynomials (which can be become more ill-conditioned in high degree cases) and instead directly compute the roots through polynomial eigenvalue problems, see, e.g., \cite{Nakatsukasa2015,doi:10.1137/15M1022513}.
\item Finally, we mention the use of quadrature compression methods, see, e.g., \cite{doi:10.1137/16M1085206}. These methods can be quite effective at significantly reducing the number of quadrature nodes of a pre-computed scheme while retaining its accuracy. However, the cost of compression can be comparatively high, hence these methods are best suited to applications in which this additional cost can be amortised.
\end{itemize}

\subsection*{Acknowledgements} This research was supported in part by the Applied Mathematics Program of the U.S.~Department of Energy (DOE) Office of Advanced Scientific Computing Research under contract number DE-AC02-05CH11231, and by a DOE Office of Science Early Career Research Program award. Some computations used resources of the National Energy Research Scientific Computing Center (NERSC), a U.S.~Department of Energy Office of Science User Facility at Lawrence Berkeley National Laboratory, operated under Contract No.~DE-AC02-05CH11231.

\bibliographystyle{siamplain}
\bibliography{references}

\end{document}

\setlength{\textfloatsep}{1.5em}

\pagestyle{myheadings}
\thispagestyle{plain}
\markboth{HIGH-ORDER QUADRATURE ON IMPLICITLY DEFINED MULTI-COMPONENT DOMAINS}{HIGH-ORDER QUADRATURE ON IMPLICITLY DEFINED MULTI-COMPONENT DOMAINS}

\maketitle

\vspace{-1em}
\centerline{\footnotesize Mathematics Group, Lawrence Berkeley National Laboratory, Berkeley, CA 94720 (\texttt{rsaye@lbl.gov})}

\vspace{2em}

\noindent This supplementary material accompanies the article \textit{High-Order Quadrature on Multi-Component Domains Implicitly Defined by Multivariate Polynomials}, published in the Journal of Computational Physics.

\vspace{1.5em}%
\centerline{\rule{0.6\linewidth}{0.4pt}}%
\vspace{1.5em}%

\section{Additional Numerical Experiments}

\subsection{Surfaces with singularities}

With the exception of the bilinear example using $\epsilon = 0$, the test problems presented in the main article have each involved polynomial functions whose zero isosurfaces give rise to a (locally) orientable interface---i.e., if $\phi$ is the polynomial which implicitly defines the interface $\Gamma = \{\phi = 0\}$, then $\nabla \phi$ is nonzero everywhere on $\Gamma$.\footnote{Technically, the randomly generated polynomials considered in \S4.3 could result in interfaces with singularities, but the set of such polynomials has measure zero and no instances were likely generated.} In the following, we examine quadrature efficacy in the case of implicitly defined geometry containing singularities, in the sense that $\nabla \phi = 0$ somewhere on $\Gamma$. %
Several examples follow. In each case, we examine the error in computing the following three integral quantities,
\[ I_{\Omega} := \int_{U \cap \{\phi < 0\}} f, \quad I_{\Gamma} := \int_{U \cap \{\phi = 0\}} f, \quad I_{\Gamma_{\textbf n}} := \int_{U \cap \{\phi = 0\}} f \mathbf n, \]
where $U$ is a $d$-dimensional cube specific to each example, and the integrand  $f : \R^d \to \R$ is given by $f(x) = \cos (\tfrac14 \|x\|_2^2)$. The quadrature algorithm is executed\footnote{We have omitted here some technical details regarding resultant calculations of polynomials having common factors or shared roots at infinity. For the purposes of quadrature, such factors can be removed to the extent they leave the interfacial geometry unaltered. In our implementation, this is handled by a variant of the main quadrature engine in which Bernstein degree reduction and GCD factoring algorithms are applied to the input and output polynomials of resultant calculations.} using the automated strategy for choosing the height direction, while the choice of one-dimensional quadrature schemes follows the guidelines derived in \S4.3.3, and shall by indicated. In each example, a $q$-refinement study is performed, using a reference solution computed with $q = 100$.

\subsubsection{Two-dimensional examples}
\label{sec:singular2d}

\begin{figure}%
\centering\includegraphics{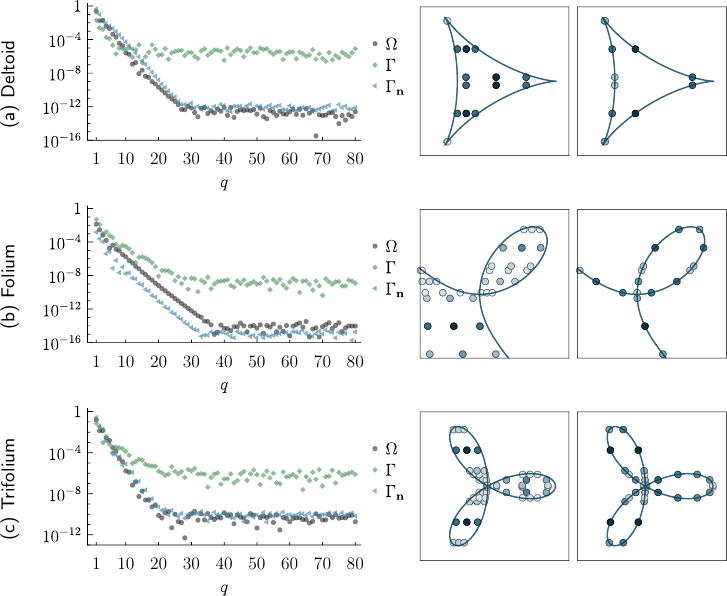}%
\caption{\textit{(left)} Quadrature accuracy as a function of $q$ for the interfacial singularity test problems of \cref{sec:singular2d}. Here, $\bullet$ denotes the relative error in the computed volumetric integral, $I_\Omega$; $\squaresymbol$ denotes the relative error in the computed surface integral in non-flux form, $I_\Gamma$; and $\trianglesymbol$ denotes the relative error in the computed flux-form surface integral, $I_{\Gamma_{\textbf n}}$. Results computed using native double precision. \textit{(right pair)}  Illustrations of the $q = 3$ volumetric (left) and surface (right) quadrature schemes; in each instance, the quadrature points are shaded according to their weight, such that low-valued weights are pale and higher-valued weights are darker.}%
\label{fig:singular2d}%
\end{figure}

\cref{fig:singular2d} illustrate three examples, defined as follows; in each example, $(\mathsf{TS},\mathsf{GL})$ quadrature is applied, i.e., tanh-sinh on the base integral, and Gauss-Legendre on the vertical line integrals when computing $I_\Omega$.

\begin{enumerate}[label=(\alph*)]
\item \textit{Deltoid.} Here, $\phi$ is given by the polynomial
\[ \phi(x,y) = (x^2 + y^2)^2 + 18 (x^2 + y^2) - 8 (x^3 - 3 xy^2) - 27, \]
and $U = (-2.5, 3.5) \times (-3, 3)$. The corresponding interface contains three cusps. According to the results shown in \cref{fig:singular2d}a, the volumetric and flux-form surface integrals each exhibit a fast convergence before plateauing around machine precision. However, the surface integral in non-flux form, $I_{\Gamma}$, quickly plateaus and does not reach the same level of accuracy. %

\vspace{0.4em} To investigate this further, numerical experiments using double-double and quad\-ruple-double precision indicate that the volumetric integral errors plateau around ${\mathcal O}(\epsilon)$, as may be expected, where $\epsilon$ is the unit round-off error. The same experiments indicate the plateau point for surface integrals is ${\mathcal O}(\epsilon^{1/3})$, which is about $10^{-5}$ in native double precision, as witnessed in \cref{fig:singular2d}a.

\vspace{0.4em} These results support the intuition that surface integrals involving cusps are highly sensitive to perturbations. Indeed, a small truncation of the cusp tip would impact a surface integral far more significantly than it would impact a volumetric integral. It is hypothesised that the ${\mathcal O}(\epsilon^{1/3})$ plateau point for surface integrals is related to the increasing root condition numbers towards the tip of the cusp, where there is a double root. Near these points, the curved surface Jacobian factors, being dependent on calculations of the surface normal, become highly sensitive to roundoff errors in the computed roots. According to our results, the volumetric integral is not sensitive to these issues.

\item \textit{Folium of Descartes.} Here, $\phi$ is given by the polynomial
\[ \phi(x,y) = x^3 + y^3 - 3 x y, \]
and $U = (-1.4, 2.1) \times (-1.5, 2)$. The corresponding interface is self-intersecting at the origin $x=y=0$. According to the results in \cref{fig:singular2d}b, the volumetric and surface-flux integrals exhibit a well-behaved, approximately exponential convergence, before plateauing around a mild multiple of machine epsilon. Meanwhile, experiments using double-double precision support the assessment that the non-flux form surface integral exhibits somewhat slower exponential convergence, with a plateau point around ${\mathcal O}(\epsilon^{1/2})$, where $\epsilon$ is the unit-roundoff error. Consequently, we observe that self-intersecting surfaces can result in a reduced maximal precision for surface integrals, albeit less severe in this example than compared to the prior cusp example of \cref{fig:singular2d}a. %

\item \textit{Trifolium.} Here, $\phi$ is given by the polynomial
\[ \phi(x,y) = (x^2 + y^2)^2 - x^3 + 3 x y^2, \]
and $U = (-1, 1.2) \times (-1.1, 1.1)$. The corresponding interface has a high degree of self-intersection at the origin $x=y=0$, see \cref{fig:singular2d}c. Once more, we observe approximately exponential convergence in $q$; however, in this example the plateau points are not as clear. According to the results shown in \cref{fig:singular2d}c, as well as experiments conducted using double-double precision, it appears that the volumetric integrals plateau around ${\mathcal O}(\epsilon^{2/3})$ and the non-flux surface integral plateaus around ${\mathcal O}(\epsilon^{1/3})$. As such, it appears the volumetric integrals are impacted by the more ill-conditioned root finding problems occurring near the junction point at the origin; indeed, the automatically chosen height direction (in the volumetric case) is $\hat x$, and for $y = 0$ the polynomial $\phi$ reduces to $\phi(x,0) = x^3(x - 1)$ where a triple root is encountered.

\end{enumerate}

\subsubsection{Three-dimensional examples}
\label{sec:singular3d}

\begin{figure}%
\centering\includegraphics{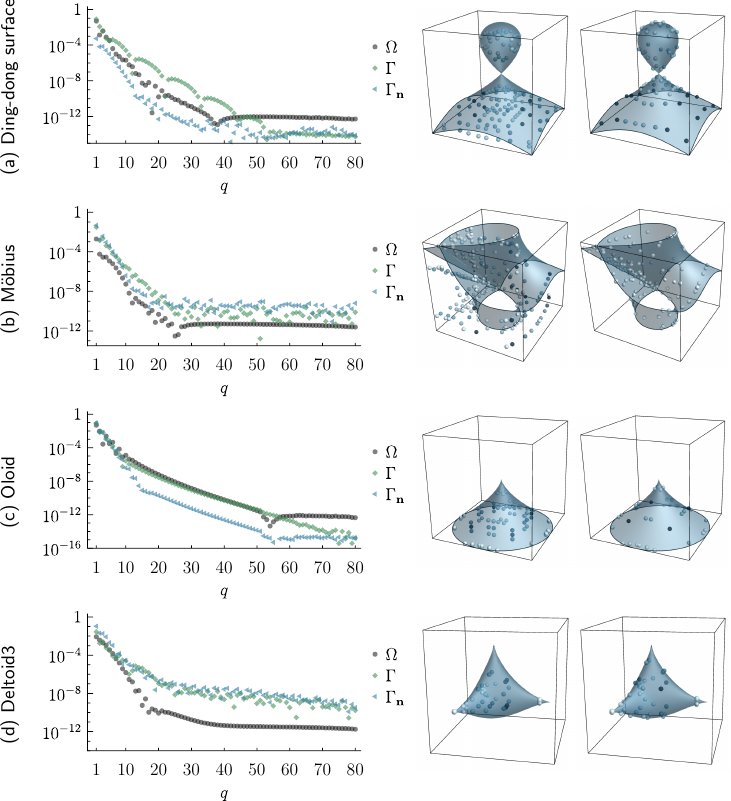}%
\caption{\textit{(left)} Quadrature accuracy as a function of $q$ for the interfacial singularity test problems of \cref{sec:singular3d}. Here, $\bullet$ denotes the relative error in the computed volumetric integral, $I_\Omega$; $\squaresymbol$ denotes the relative error in the computed surface integral in non-flux form, $I_\Gamma$; and $\trianglesymbol$ denotes the relative error in the computed flux-form surface integral, $I_{\Gamma_{\textbf n}}$. Results computed using native double precision. \textit{(right pair)}  Illustrations of the volumetric (left) and surface (right) quadrature schemes, with $q = 2$, $2$, $3$, and $2$, respectively from top to bottom; in each instance, the quadrature points are shaded according to their weight, such that low-valued weights are pale and higher-valued weights are darker.}%
\label{fig:singular3d}%
\end{figure}

We consider four examples in 3D, as illustrated in \cref{fig:singular3d}. In each example, $(\mathsf{TS},\mathsf{TS},\mathsf{GL})$ quadrature is applied, i.e., tanh-sinh on every base integral, and Gauss-Legendre on the vertical line integrals over the first eliminated axis when computing $I_\Omega$. Owing to the degree of complexity, the analysis of finite precision effects is less clear for these examples compared to the 2D examples. As such, only rough conclusions are drawn.

\begin{enumerate}[label=(\alph*)]
\item \textit{Ding-dong surface}, a cubic surface of revolution. Here, $\phi$ is given by the polynomial
\[ \phi(x,y,z) = x^2 + y^2 - (1-z)z^2 \]
and $U = (-1,1)^3$. The corresponding interface has a single point of singularity at the origin $x = y = z = 0$. According to the results in \cref{fig:singular3d}a, approximately exponential convergence rates are obtained, for all three integral quantities, with the non-flux surface integral showing the slowest rate of convergence.

\item \textit{M\"obius surface}, a portion of which yields the well-known M\"obius strip. Here, $\phi$ is given by the polynomial 
\[ \phi(x,y,z) = -4 y + x^2 y + y^3 + 4 x z - 2 x^2 z - 2 y^2 z + y z^2 \]
and $U = (-1.25, 2.75) \times (-2, 2) \times (-2,2)$. The corresponding interface has a line of self-intersection. The results shown in \cref{fig:singular3d}b indicate approximately exponential convergence rates, up to the plateau points shown. %

\item \textit{Oloid.} Here, $\phi$ is given by the polynomial 
\[ \phi(x,y,z) = x^2 + y^2 + z^3 \]
and $U = (-1,1)^3$. The corresponding interface is a surface of revolution exhibiting a single point of singularity, in fact, a cusp, at the origin $x = y = z = 0$. According to the results in \cref{fig:singular3d}c, this cusp does not pose a significant problem, with approximately exponential asymptotic convergence rates obtained. %

\item \textit{Deltoid3.} Here, $\phi$ is given by the polynomial
\[ \phi(x,y,z) = (x^2 + y^2 + z^2)^2 + 18 (x^2 + y^2 + z^2) - 8 (x^3 + z^3 - 3 x y^2) - 27 \]
and $U = (-2.5, 3.5) \times (-3, 3) \times (-2, 4)$. The corresponding interface is a kind of 3D generalisation of the 2D deltoid, hence the name given to it here; it exhibits four cusps. This surface is perhaps the most challenging example presented in this work, and stress-tests many aspects of the quadrature algorithm.

\end{enumerate}

\subsubsection{Addendum}

We emphasise that, in all of the above 2D and 3D test problems involving surface singularities, the black-box nature of the quadrature algorithm means that it is in some ways ``blind'' to the geometry at hand. In particular, no specific attention is given to the ill-conditioned nature of root finding near self-intersections, cusps, and other points of singularity; moreover, no specific attention is given to the kinds of interfacial curvature encountered along the way. All of the above results could likely be improved if the implementation is tailored to the various intricacies at play. Our main purpose for these test problems is to stress-test the quadrature algorithm and demonstrate its robustness in these settings, yielding, at a minimum, modest accuracy for small quadrature orders $q$.

\subsection{Subdivision of randomly generated geometry}

In this section, we explore the effectiveness of applying a small amount of constraint cell subdivision on the randomly generated geometry considered in \S4.3 of the main article. As discussed therein, the quadrature algorithm is designed so as to not rely on this kind of subdivision; nevertheless, a small amount of subdivision permits a more locally-adaptive choice of height direction/elimination axis, and can result in an increased likelihood of rapid convergence via Gauss-Legendre quadrature.

We repeat here the experiments of \S4.3, but with two modifications: (i) instead of applying the quadrature algorithm directly to $U = (-1,1)^d$, we apply it to each subcell arising from a uniform Cartesian subdivision of $U$ into $k \times \cdots \times k$ subcells, accumulating the results so as to compute the mass matrix $I^{\pm}$ on the whole of $U$; and (ii) on each subcell, the quadrature scheme is chosen automatically, such that Gauss-Legendre is utilised on the base integral whenever it can be proven there are no branching points/vertical tangents, and tanh-sinh is utilised in all other cases. It follows that some subcells will utilise different height directions compared to other subcells. As $k$ is made larger, the interfacial geometry on each subcell appears locally flatter, and the proportion of subcells in which Gauss-Legendre is effective will increase.

\begin{figure}%
\centering\includegraphics{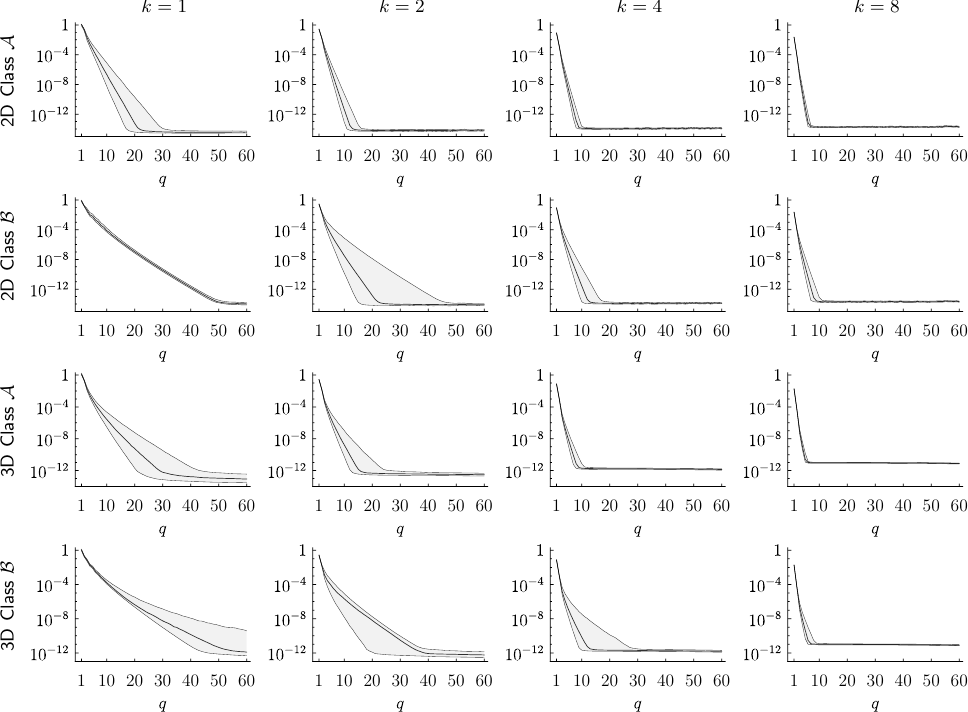}%
\caption{Quadrature accuracy as a function of $q$ and subdivision amount $k$, for the randomly generated geometry considered in \S4.3 of the main article; recall that instances of class $\mathcal A$ typically involve nicer interfacial geometry, whereas those of $\mathcal B$ have complex, high curvature pieces. Each column corresponds to a different subdivision amount, and each row corresponds to aggregated statistics of 1000 randomly generated polynomials (generated once and then held fixed across each row). The shaded region demarcates the boundaries of the first and third quartiles of the error metric $E_q$ (measuring the relative error in computed mass matrices), while the interior line indicates the median.}
\label{fig:subdivision}
\end{figure}

\cref{fig:subdivision} displays the results of this study, for $k \in \{1,2,4,8\}$ (the results for $k = 1$ are essentially identical to those of the main article and are included for comparison). As expected, increasing the subdivision amount results in faster convergence, across all cases in 2D and 3D. For class $\mathcal A$ randomly generated geometry, which typically involves smooth interface pieces, increasing $k$ results in a monotonic improvement, both in the median as well as a monotonic reduction of the interquartile spread. For class $\mathcal B$ polynomials, which typically involve complex, high curvature pieces, increasing $k$ from $1$ to $2$ reduces the spread, but not significantly; this is because a subdivision by half in each axis is typically not enough for class $\mathcal B$ geometry, in the sense that some subcells continue to have branching points in every possible height direction (thus requiring tanh-sinh quadrature). With $k = 4$, there is a better chance that every subcell can find a height direction in which Gauss-Legendre is appropriate, and we see that with $k = 8$, the spread in the 3D class $\mathcal B$ results has significantly reduced. %

In summary, it can be seen that a small amount of subdivision can be effective in substantially reducing the integral error for the kinds of randomly generated interfacial geometry depicted in Fig.~6 of the main article. There are many possibilities here to further automate this process via other kinds of adaptive methods.

\section{Masking Algorithms} %
\label{app:orthant}

As motivated in \S3 of the main article, a number of operations in the quadrature algorithm benefit from a cheap and robust means of isolating zeros of a  multivariate polynomial and to detect when two polynomials have non-intersecting zero isosurfaces. A simple and effective approach achieving these objectives is through ``orthant tests'' combined with a small amount of subdivision to exploit the effectiveness of Bernstein polynomial range evaluation; we discuss these algorithms here.

\subsection{Orthant tests} Suppose that $f(x) = \smash{\sum_i f_i \prod_{\ell=1}^d b_{i_\ell}^{n_\ell} (x_\ell)}$ and $g(x) = \smash{\sum_i g_i \prod_{\ell=1}^d b_{i_\ell}^{n_\ell} (x_\ell)}$ are two Bernstein polynomials of equal degree $(n_1,\ldots,n_d)$. If we can prove that a linear combination of $f$ and $g$ is nonzero on $[0,1]^d$, then it follows that $f$ and $g$ do not share a root in that domain. As a special case, let us determine whether we can (cheaply) find a constant $\alpha \in \R$ such that $f + \alpha g > 0$ on $[0,1]^d$. Using the convex hull/variation diminishing property of Bernstein polynomials, we can do so by finding, if possible, $\alpha \in \R$ such that $f_i + \alpha g_i > 0$ for all multi-indices $i$. We call this an \textit{orthant test} since it is asking if the vector of coefficients of $f + \alpha g$ is in the positive orthant. The test is a particularly simple case of a linear programming problem and a fast algorithm for determining if $\alpha$ exists is simple to design: first, let $\alpha = (-\infty,\infty)$; then, for each $i$, progressively restrict $\alpha$ to the range of values for which $f_i + \alpha g_i > 0$; upon conclusion, if the interval is non-empty then the orthant test ``succeeded'' since we know there exists an $\alpha$ such that $f + \alpha g > 0$. Swapping the role of $f$ and $g$ then leads to a simple algorithm to determine if there exists $\alpha, \beta \in \R$ such that $\alpha f + \beta g > 0$. It is important to note these tests are designed to cheaply determine if $f$ and $g$ do not share a common root; a negative result (i.e., no such $\alpha$ or $\beta$ satisfies the orthant test) does not imply that a common root exists. However, these tests are closely related to the effectiveness of Bernstein polynomials in evaluating polynomial ranges or bounds; in particular, they become more and more accurate under degree elevation and subdivision operations.

\subsection{Mixed degree orthant tests} In the above, we assumed that $f$ and $g$ had the same Bernstein degree. In the case that they do not, it is simple to elevate their degrees (in each dimension) so that they match, and then apply the equal-degree orthant test. (Of note is that degree elevation actually increases the accuracy of Bernstein range evaluation \cite{garloff1993bernstein,lin1995methods}.) If the test succeeds, then the polynomials $f$ and $g$ do not share a common root in $[0,1]^d$; if the result is negative, then $f$ and $g$ may or may not share a common root.

\subsection{Subdivision} The above orthant tests can benefit from a small amount of subdivision to increase their certitude, i.e., decrease the number of cases in which the orthant test yields an uncertain result. We accomplish this by dividing the unit cube $[0,1]^d$ into a regular grid of $\mathcal M \times \cdots \times \mathcal M$ subcells, and then, for each subcell, apply the orthant test with input polynomials $f$ and $g$ computed relative to that subcell (using the de Casteljau algorithm). In particular, when $\mathcal M$ is a power of two, a recursive subdivision strategy can be applied for improved efficiency, as follows. If the orthant test applied to $f$ and $g$ proves there are no shared roots in $[0,1]^d$, then no subdivision is needed and an empty mask can be returned. Otherwise, divide the current cube into $2^d$ half-cubes and apply the orthant test to each one: whenever the orthant test succeeds, the subdivision process for that piece can terminate; else, the subdivision process recursively continues, up to the maximum level corresponding to $\mathcal M \times \cdots \times \mathcal M$ subcells.

\begin{algorithm}[t]%
\caption{\small $\texttt{intersectionMask}\bigl((f,f_m),(g,g_m)\bigr)$. In $d$ dimensions, given two Bernstein polynomials $f(x) = \smash{\sum_i f_i \prod_{\ell=1}^d b_{i_\ell}^{n_\ell} (x_\ell)}$ and $g(x) = \smash{\sum_i g_i \prod_{\ell=1}^d b_{i_\ell}^{m_\ell} (x_\ell)}$, and their masks $f_m, g_m : {\mathbb N}_{\mathcal M}^d \to \{0,1\}$, compute a mask analysing the intersection of  their zero isosurfaces.}%
\label{algo:intersectionMask}%
\begin{algorithmic}[1]%
	\State Initialise the mask $m :  {\mathbb N}_{\mathcal M}^d \to \{0,1\}$ as empty, $m \equiv 0$.
	\State Execute recursive process by calling the following procedure with $a = 0$ and $b = ({\mathcal M},\ldots,{\mathcal M})$.
	\Procedure{examineRegion}{$a \in {\mathbb N}^d$, $b \in {\mathbb N}^d$}
		\Statex \textit{\hspace{0.5em}Test whether $f$ and $g$ have overlapping masks in the subrectangle $[a,b]$:}
		\State \textbf{if} $f_m(i) g_m(i) = 0$ for all $a < i \leq b$ \textbf{then return}
		\Statex \textit{\hspace{0.5em}Apply orthant test in the subrectangle corresponding to $[a,b]$:}
		\State Let $x_a := a / {\mathcal M} - \epsilon$ and $x_b := b / {\mathcal M} + \epsilon$.
		\parState{Apply a multi-dimensional de Casteljau algorithm to compute the Bernstein coefficients of $f$ and $g$ relative to the subrectangle $[x_a,x_b]$, denoting the result $f_{ab}$ and $g_{ab}$.}
		\parState{If $f_{ab}$ and $g_{ab}$ are not of identical degree, apply a multi-dimensional Bernstein degree elevation algorithm so that they match.}
		\parState{Apply orthant test to determine whether there exists  $\alpha,\beta \in \R$ such that $\alpha f_{ab,i} + \beta g_{ab,i} > 0$ for every coefficient $i$.}
		\State \textbf{if} orthant test succeeded \textbf{then return}
		\Statex \textit{\hspace{0.5em}If at terminus of recursive subdivision process, set the mask to true and return:}
		\If{$b = a + 1$}
			\State Set $m(b) := 1$
			\State \textbf{return}
		\EndIf
		\Statex \textit{\hspace{0.5em}Otherwise, subdivide $[a,b]$ into $2^n$ equal subrectangles and apply recursion:}
		\State Set $\delta := (b - a)/2$.
		\For{every $i \in \{0,1\}^d$}
			\State \textbf{call} \textsc{examineRegion}($a + i\,\delta$, $a + (i + 1)  \delta$).
		\EndFor
	\EndProcedure
	\State \textbf{return} $m$.
	\end{algorithmic}%
\end{algorithm}

\subsection{Algorithm} The above ideas, using orthant tests, degree elevation when necessary, and recursive subdivision, are combined to yield the algorithm outlined in \cref{algo:intersectionMask}, implementing the routine $\texttt{intersectionMask}\bigl((f,f_m),(g,g_m)\bigr)$ discussed in \S3.3. A special case of this algorithm, with $g \equiv 0$ and $g_m \equiv 1$, can be used to implement $\verb|nonzeroMask|(f,f_m)$. In this work, the overlap threshold, which is used to eliminate edge cases that could be affected by roundoff error (e.g., a polynomial has a root exactly on the border of two adjacent subcells), is set via $\epsilon = 2^{-6} {\mathcal M}^{-1}$, i.e., about 1\% of the width of a subcell. This tolerance is relatively lenient, but otherwise not of any particular significance.

\subsection{Discussion} These masking operations are designed to be robust, simple, and cheap. A number of other possibilities could be devised to examine whether a pair of multivariate polynomials share common real roots in the unit cube $[0,1]^d$. As examples:
\begin{itemize}
\item One approach may be to build two masks, one mask that determines the subcells for which $f$ is nonzero, the second mask performing the same role for $g$, and then binary-\textsc{and} these masks to determine the subcells in which $f$ and $g$ could share a common root. However, this approach is nowhere near as accurate as the orthant test approach discussed above. As a trivial example, $f(x) = x$ and $g(x) = x + 0.001$ do not share a common root, yet their masks would both be nonempty in the subcell closest to the origin (unless $\mathcal M$ was excessively large). The orthant test, on the other hand, is free to build \textit{any} linear combination of $f$ and $g$, whose Bernstein coefficients are uniform in sign, to conclude no shared roots exist. %

\item As another approach, the idea of the orthant test could be generalised to determine whether there exists any strictly-positive polynomials $\alpha$ and $\beta$ such that $\alpha f + \beta g > 0$ on $[0,1]^d$. (This idea is closely related to positivstellensatz methods in real algebraic geometry.) By examining the Bernstein coefficients of $\alpha f + \beta g$ (in a suitable basis of sufficient degree), this approach would then lead to a linear programming problem in the coefficients of $\alpha$ and $\beta$. This idea was explored as part of this work: while the accuracy of these methods is indeed higher (e.g., in producing fewer false assertions that shared roots probably exist), it was found that the complexity and expense incurred from setting up and solving the linear program did not justify the adoption of this method. Instead, the simple orthant test, in which $\alpha$ and $\beta$ are constants, yields a robust and effective means of masking.
\end{itemize}

\section{Resultant Computations}
\label{app:resultant}

\begin{algorithm}[t]%
\caption{\small $R(f,g;k)$. In $d$ dimensions, given two Bernstein polynomials $f(x) = \smash{\sum_i f_i \prod_{\ell=1}^d b_{i_\ell}^{n_\ell} (x_\ell)}$ and $g(x) = \smash{\sum_i g_i \prod_{\ell=1}^d b_{i_\ell}^{m_\ell} (x_\ell)}$, and an elimination axis $k$, compute the resultant of $f$ and $g$ in the direction $k$.}%
\label{algo:resultant}%
\begin{algorithmic}[1]
	\State Let $r = (r_1, \ldots, r_{k-1}, r_{k+1}, \ldots, r_d)$ denote the maximum possible degree of the resultant (see eq.\,(6), main article).
	\For{every multi-index $i \in \prod_{\ell = 1, \ell \neq k}^d \{0, \ldots, r_\ell\}$}
		\parState{Let ${\bar x}_i \in [0,1]^{d-1}$ be the $i$th Chebyshev tensor-product grid point inclusive of endpoints, i.e., ${\bar x}_{i_\ell} = \tfrac12 + \tfrac12 \cos \tfrac{i_\ell}{r_\ell} \pi$.}
		\parState{Let $f_{\bar x}$ and $g_{\bar x}$ be the line restrictions of $f$ and $g$ at $\bar x$, i.e., the degree $n_k$ and $m_k$ Bernstein polynomials such that $f_{\bar x}(y) = f(\bar x + y e_k)$ and $g_{\bar x}(y) = g(\bar x + y e_k)$.}
		\If{$n_k = m_k$}
			\State Compute the B\'ezout matrix $A := B(f_{\bar x}, g_{\bar x})$.
		\Else
			\State Compute the Sylvester matrix $A := S(f_{\bar x}, g_{\bar x})$.
		\EndIf
		\State Define $q_i := \det A$, computed via QR with column pivoting.
	\EndFor
	\State \textbf{return} the unique Bernstein polynomial of degree $r$ which interpolates $q_i$ at the nodes ${\bar x}_i$.
	\end{algorithmic}%
\end{algorithm}

In this work, the resultant of two multivariate polynomials is computed by interpolating the nodal values of a one-dimensional resultant. In particular, applying the relation
\[ \bigl[R(f, g; k)\bigr] (\bar x) = R(f_{\bar x}, g_{\bar x}) \text{ for all } \bar x, \]
we evaluate the right hand side at the nodes of a tensor-product grid, and use these values to recover the polynomial $R(f,g;k)$ in Bernstein form. \cref{algo:resultant} implements this algorithm with the following approach:
\begin{itemize}
\item To evaluate the one-dimensional resultants $R(f_{\bar x}, g_{\bar x})$ at a given base point $\bar x$, the univariate polynomials $f_{\bar x}$ and $g_{\bar x}$ are found by a simple line restriction of the given multivariate polynomials $f$ and $g$; a standard QR algorithm with column pivoting is then used to compute the determinant of the Sylvester or B\'ezout matrix of $f_{\bar x}$ and $g_{\bar x}$.
\item There is some freedom in choosing the interpolation points; however, as is often the case in polynomial recovery problems, it is best to use a non-uniform distribution. Here, a Chebyshev scheme, inclusive of the domain endpoints, has been used. In particular, the standard Chebyshev points, relative to the interval $[-1,1]$, are given by $\{\cos \frac{2\ell - 1}{2n} \pi\}_{\ell=1}^n$ and are non-inclusive of the endpoints; we instead use the points $\{\cos \frac{\ell - 1}{n - 1} \pi\}_{\ell=1}^n$. The numerical conditioning of the associated Vandermonde matrix is slightly improved for the endpoint-inclusive scheme, though both Chebyshev schemes have the same asymptotic conditioning as the number of nodal points increases.
\item After evaluation of the 1D resultants on a tensor-product grid, the task is then to compute a Bernstein polynomial interpolating these values. One approach to Bernstein polynomial recovery is through the simple-to-implement and fast algorithm developed by Ainsworth and S\'anchez \cite{doi:10.1137/15M1046113}; however, we observed strong instabilities in this method for high degree resultant calculations. Instead, we have adopted here a multi-dimensional SVD method. The method exploits the tensor-product formulation to recast the Bernstein polynomial recovery problem as a sequence of one-dimensional Bernstein polynomial recovery problems; if the degree of multivariate polynomial is $(n_1,\ldots,n_d)$, then the first 1D problem is of degree $n_1$, the second is of degree $n_2$, etc. In each one-dimensional problem, the associated Vandermonde matrix $V = (b_{j-1}^{n}(x_i))_{i,j=1}^{n+1}$ is formed, its SVD computed, and then used to solve the Vandermonde system. Importantly, for each degree $n$ one can precompute or cache the SVD for multiple interpolation operations, since the node distribution for a particular $n$ does not change.
\item Numerical experiments indicate that the condition number of the tensor-product Vandermonde system in $d$ dimensions, corresponding to the recovery of a degree $(n_1,\ldots,n_d)$ Bernstein polynomial, is approximately asymptotically $\prod_{\ell = 1}^d 2^{n_\ell}$. Clearly, this conditioning could be problematic for high-degree polynomials. In our implementation a simple regularisation is used whereby the SVD is truncated and a pseudoinverse applied, for any singular values below a modest multiple of machine epsilon. By design, the quadrature algorithms developed in this work are quite robust to this kind of regularisation, and our experiments indicate that the quadrature accuracy is not strongly affected by this regularisation, at least for the class of polynomials considered in this work. %
\end{itemize}

\section{Computing Roots of Univariate Bernstein Polynomials}
\label{app:roots}

A variety of approaches have been developed to compute the roots of univariate Bernstein polynomials, including through computation of companion matrices and their eigenvalues \cite{WINKLER2000179,WINKLER2003153}; using Bernstein subdivision algorithms and convex hull properties to first isolate (real) roots and then polish with an iterative method such as Newton's \cite{Farouki2008,Spencer}; and through generalized eigenvalue problems \cite{doi:10.1137/S0895479899365720,10.1007/978-3-030-41258-6_6}, among a number of other possibilities. For a discussion of some of these methods, see the dissertation of Spencer \cite{Spencer}. In this work, we have adopted a hybrid approach, using a combination of Bernstein subdivision and Newton's method (for speed), falling back on a backward-stable generalised eigenvalue method, as follows.

Let $f$ be a given univariate Bernstein polynomial of degree $n$, where $f(x) = \sum_{i=0}^n f_i b_i^n(x)$. Two key properties of Bernstein polynomials concerning roots are as follows---assuming none of the coefficients $f_i$ are zero,\footnote{Zero-valued coefficients are an edge case in Descartes' rule of signs and require extra care.}
\begin{enumerate}[label=(\roman*)]
\item if every coefficient $f_i$ has the same sign, then there are no real roots of $f$ in $[0,1]$; 
\item if there is only one sign change, meaning there exists $\ell$ such that $\text{sign}(f_i) \neq \text{sign}(f_j)$ for all $0 \leq i \leq \ell < j \leq n$, then there is exactly one (simple) real root of $f$ in $(0,1)$.
\end{enumerate}
These properties naturally suggest a simple recursive algorithm for computing roots. Beginning with the interval $(0,1)$, check whether either condition holds; if not, subdivide the interval into two halves (computing the Bernstein coefficients of $f$ relative to each half using the de Casteljau algorithm), test these conditions on each half, and continue recursively if necessary up to a maximum subdivision level. On any subinterval satisfying property (ii), we can then apply Newton's method to compute the associated root. If Newton's method fails to converge in a limited number of iterations, or if the subdivision process reaches a maximum level of recursion, then we fall-back to an eigenvalue method developed by J\'onsson \cite{doi:10.1137/S0895479899365720}, as follows.

The generalised eigenvalue method of J\'onsson reformulates the Bernstein root finding problem as follows: $\lambda$ is a root of $f$ if and only if $A - \lambda B$ is singular, where the $n \times n$ matrices $A$ and $B$ are given by
\[ A = \begin{bmatrix} 0 & 1 \\ & 0 & 1 \\ & & \ddots & \ddots \\ & & & 0 & 1 \\ -f_0 & -f_1 & \cdots & -f_{n-2} & -f_{n-1} \end{bmatrix}, \]
and
\[ B = \begin{bmatrix} \tfrac{n}{1} & 1 \\ & \tfrac{n-1}{2} & 1 \\ & & \ddots & \ddots \\ & & & \tfrac{2}{n-1} & \hspace{-3em}1 \\ -f_0 & -f_1 & \cdots & -f_{n-2} & -f_{n-1} + \tfrac1n f_n \end{bmatrix}. \]
The complete set of roots of $f$ (real or complex) can then be computed using a generalised eigenvalue method applied to $A x= \lambda B x$, such as the \textsc{Lapack} routine \texttt{dggevx}. As discussed in \cite{doi:10.1137/S0895479899365720}, this method of computing roots of Bernstein polynomials is particularly well-conditioned, even for tiny coefficient values, and is backward stable in the sense that the computed roots are the exact roots of a nearby polynomial with nearby coefficients. For improved accuracy, we note that it is sometimes beneficial to apply balancing transforms (as provided, e.g., by \textsc{Lapack}) to gain further improvements in conditioning.

To summarise and provide a few more implementation details, the root finding algorithms used in this work operate as follows:
\begin{itemize}
\item Given a univariate polynomial $f$ of Bernstein degree $n > 2$, we first try for a fast method. Using the de Casteljau algorithm, the method recursively subdivides the interval $(0,1)$ into subintervals, until either: (a) it can be proven, using solely examination of the signs of the Bernstein coefficients, whether there is exactly 0 or 1 real roots on each subinterval, or (b) a user-defined maximal recursion depth is met.
\item On each subinterval containing exactly one real root, we apply Newton's method, safeguarded by bisection. Numerical experiments indicate that, in the vast majority of cases, Newton's method finds the root with $\approx 16$ digits of accuracy in at most 8 iterations, and $\approx 64$ digits of accuracy in at most 10 iterations. Correspondingly, we implemented a hard limit of 12 iterations in Newton's method. If the method fails to converge to the root by this time, then we fall-back to the eigenvalue method (though this very rarely occurs).
\item The fast subdivision method may not always work, e.g., in the case of (near) double roots, or in high degree cases. Its associated maximum recursion depth is a user-defined parameter controlling the tradeoff between excessive subdivision and application of the eigenvalue method. We have heuristically set this limit to be a maximum of 4 levels of subdivision, i.e., the smallest subinterval is no smaller than a 16th of the original size, though more adaptive strategies could be developed. %
\item In rare circumstances, the subdivision method may yield Bernstein polynomials with one or more coefficients very close to zero; this situation could arise, e.g., when a root is very near the endpoint of a subinterval. In this circumstance, it is ill-advised to apply Descartes' rule of signs to check conditions (i) and (ii) above. A tolerance is used to check for such cases, and if close-to-zero coefficients are detected, we fall-back to the eigenvalue method.
\item When the eigenvalue method is invoked, we use it to compute all roots of the polynomial, both real and complex. A tolerance is applied to filter out all non-real roots, and then only those in the interval $(0,1)$ are returned by the method.
\item Finally, in the simple cases of input Bernstein polynomials of degree $n = 1$ or $n = 2$, we bypass the above methods and directly apply standard routines for computing the roots of linear or quadratic polynomials.
\item Overall, we found this hybrid approach to be quite effective for the typical polynomials produced the quadrature algorithms. In a majority of cases, the fast subdivision-based approach using Newton's method is successful at finding (and computing to full machine precision) all real roots in $(0,1)$. Less often, and usually only for higher degree polynomials, is the more expensive eigenvalue method executed. We also found these methods to be insensitive to ``roots at infinity'', i.e., when the input polynomial has actual degree less than the degree used in its Bernstein representation.
\end{itemize}

\section{Tanh-sinh Quadrature}
\label{app:tanhsinh}

A tanh-sinh quadrature scheme, relative to the interval $(-1,1)$, using $q > 1$ points, takes the form
\[ \int_{-1}^1 f(x)\,dx \approx \sum_{\ell = 1}^q w_\ell f(x_\ell), \]
where
\[ x_\ell = \tanh \bigl( \tfrac12 \pi \sinh t_\ell \bigr), \quad w_\ell = \frac{\tfrac12 h_q \pi \cosh t_\ell}{ \cosh^2 (\tfrac12 \pi \sinh t_\ell)}, \]
are the quadrature nodes and weights. Here, $\smash{t_\ell = h_q \left\lfloor \tfrac{\ell}{2} \right\rfloor (-1)^\ell}$ when $q$ is odd, and $\smash{t_\ell = h_q (\left\lceil \tfrac{\ell}{2} \right\rceil - \tfrac12) (-1)^\ell}$ when $q$ is even, for $\ell = 1, \ldots, q$. The step size $h_q$ is associated with the tanh-sinh transformation to a trapezoidal rule on the whole real line $(-\infty,\infty)$; its value is usually chosen as a balance between the discretization errors of the trapezoidal rule and the need to truncate it to a finite number of quadrature points. We follow here the suggestions of \cite{doi:10.1137/130932132,vanherck2020tanhsinh} to define
\begin{equation} \label{eq:h} h_q = \tfrac{2}{q} W\bigl(C \pi (q - 1)\bigr) \end{equation}
where $W(z)$ is the Lambert W-function. The constant $C$ relates to the expected strip of analyticity of the transformed integrand $f$; in many applications, however, this can be difficult to characterise. Often, a good choice is to take $C = 1$ \cite{doi:10.1137/130932132,vanherck2020tanhsinh}. We found in this work that the quadrature schemes (for implicitly defined geometry) can yield better accuracy if a somewhat smaller value of $C$ is used; in particular, in \eqref{eq:h} we have set $C = 0.6$, determined empirically from specific kinds of geometry as well as randomly generated interfaces. It is hypothesised this reduction is related to the typical behaviour of the multi-valued height functions produced by the dimension reduction process, specifically the decay-to-zero endpoint singularities at the branching points/vertical tangents of the interface $\Gamma$.

For the actual application of these schemes to the quadrature algorithms for implicitly defined geometry, we slightly modify the above scheme by enforcing a normalisation of the weights. The standard tanh-sinh quadrature scheme does not integrate constants exactly, but a simple normalisation fixes this: the weights $w_\ell$ above are replaced by $\tilde{w}_\ell$, where
\[ \tilde{w}_\ell = w_\ell \frac{2}{\sum_{j=1}^q w_j}. \]
In the dimension reduction quadrature algorithm, this normalisation is used to enforce the requirement of exactly integrating constant functions. Critically, this normalisation does not impact the (near) exponential convergence of tanh-sinh schemes (when applied to sufficiently nice integrands).

We have omitted here some implementation details which can be important for integrands with blow-up endpoint singularities, see, e.g., \cite{vanherck2020tanhsinh,baileynotes}; these details relate to potential underflow in the calculation of weights, or the floating point rounding of the quadrature points closest to the interval endpoints $\pm 1$. However, in the application of this work, no integrands exhibit blowup singularities and these particularities are not of major concern.

Finally, we point out the above formulation generates quadrature schemes for $q > 1$ points. In the special case of $q = 1$, the tanh-sinh scheme used in our numerical experiments is defined to be the same as the midpoint rule.

\bibliographystyle{siamplain}
\bibliography{references}